\documentclass[letterpaper, 10 pt, conference]{ieeeconf}%
\usepackage{makeidx,ifmtarg}
\usepackage{graphicx}
\usepackage{multicol}
\usepackage[bottom]{footmisc}
\usepackage{latexsym}
\usepackage{cite,color,comment,xspace}
\usepackage{amsfonts}
\usepackage{amsmath}
\usepackage{mathrsfs}
\usepackage{amssymb}
\usepackage{algorithm}
\usepackage{algorithmic}
\usepackage[skip=2pt,font=scriptsize]{caption}
\usepackage{subcaption}
\usepackage{authblk}
\usepackage[bottom]{footmisc}
\usepackage{latexsym}
\usepackage{epstopdf}%
\providecommand{\U}[1]{\protect\rule{.1in}{.1in}}
\IEEEoverridecommandlockouts
\DeclareMathOperator*{\argmax}{\arg\!\max}
\allowdisplaybreaks
\begin{document}

\title{{\LARGE \textbf{A Submodularity-Based Approach for Multi-Agent Optimal
Coverage Problems}}}
\author{Xinmiao Sun, Christos G. Cassandras and Xiangyu Meng \thanks{This work was
supported in part by NSF under grants ECCS-1509084 and IIP-1430145, by AFOSR
under grant FA9550-12-1-0113, and by the MathWorks.}\thanks{The authors are
with the Division of Systems Engineering and Center for Information and
Systems Engineering, Boston University, Brookline, MA 02446, USA
\texttt{{\small \{xmsun,cgc,xymeng\}@bu.edu}}}}
\maketitle

\begin{abstract}
We consider the optimal coverage problem where a multi-agent network is
deployed in an environment with obstacles to maximize a joint event detection
probability. The objective function of this problem is non-convex and no
global optimum is guaranteed by gradient-based algorithms developed to date.
We first show that the objective function is monotone submodular, a class of
functions for which a simple greedy algorithm is known to be within $1-1/e$
of the optimal solution. We then derive two tighter lower bounds by exploiting
the curvature information (total curvature and elemental curvature) of the
objective function. We further show that the tightness of these lower bounds
is complementary with respect to the sensing capabilities of the agents. The
greedy algorithm solution can be subsequently used as an initial point for a
gradient-based algorithm to obtain solutions even closer to the global
optimum. Simulation results show that this approach leads to significantly
better performance relative to previously used algorithms.

\end{abstract}

\section{INTRODUCTION}

Multi-agent systems involve a team of agents, e.g., vehicles, robots, or
sensor nodes, that cooperatively perform one or more tasks in a mission space
which may contain uncertainties in the form of obstacles or random event
occurrences. Examples of such tasks include environmental monitoring,
surveillance, or animal population studies among many. Optimization problems
formulated in the context of multi-agent systems, more often than not, involve
non-convex objective functions resulting in potential local optima, while
global optimality cannot be easily guaranteed.

One of the fundamental problems in multi-agent systems is the optimal coverage
problem where agents are deployed so as to cooperatively maximize the coverage
of a given mission space
\cite{SM2011,cgc2005,caicedo2008coverage,caicedo2008performing,breitenmoser2010voronoi}
where \textquotedblleft coverage\textquotedblright\ is measured in a variety
of ways, e.g., through a joint detection probability of random events
cooperatively detected by the agents. The problem can be solved by either
on-line or off-line methods. Some widely used on-line methods, such as
distributed gradient-based algorithms~\cite{cgc2005,Minyi2011,Gusrialdi2011}
and Voronoi-partition-based
algorithms~\cite{CM2004,gusrialdi2008voronoi,breitenmoser2010voronoi},
typically result in locally optimal solutions, hence possibly poor
performance. To escape such local optima, a \textquotedblleft boosting
function\textquotedblright\ approach is proposed in~\cite{Sun2014} whose
performance can be ensured to be no less than that of these local optima.
Alternatively, a \textquotedblleft ladybug exploration\textquotedblright%
\ strategy is applied to an adaptive controller in~\cite{schwager2008}, which
aims at balancing coverage and exploration. However, these on-line approaches
cannot quantify the gap between the local optima they attain and the global
optimum. Off-line algorithms, such as simulated annealing
\cite{van1987simulated,bertsimas1993simulated}, can, under certain conditions,
converge to a global optimal solution in probability. However, they are
limited by a high computational load and slow convergence rate.

Related to the optimal coverage problem is the \textquotedblleft maximum
coverage\textquotedblright\ problem
\cite{khuller1999budgeted,berman2002generalized}, where a collection of
discrete sets is given (the sets may have some elements in common and the
number of elements is finite) and at most $N$ of these sets are selected so
that their union has maximal size (cardinality). The objective function in the
maximum coverage problem is \emph{submodular}, a special class of set
functions with attractive properties one can exploit. In particular, a well
known result in the submodularity theory \cite{nemhauser1978analysis} is the
existence of a lower bound for the global optimum provided by any feasible
solution obtained by the \textit{greedy algorithm}, i.e., an algorithm which
iteratively picks the set that covers the maximum number of uncovered elements
at each iterative step. Defining, for any integer number $N$ of sets,
$L(N)=f/f^{\star}$ where $f^{\star}$ is the global optimum and $f$ is a
feasible solution obtained by the greedy algorithm, it is shown in
\cite{nemhauser1978analysis} that $L(N)\geq1-\frac{1}{e}$. In other words,
since $f^{\star}\leq(1-\frac{1}{e})^{-1}f$, one can quantify the optimality
gap associated with a given solution $f$.

In our past work \cite{Sun2014}, we studied the optimal coverage problem with
agents allowed to be positioned at any feasible point in the mission space
(which generally includes several obstacles) and used a distributed
gradient-based algorithm to determine optimal agent locations. Depending on
initial conditions, a trajectory generated by such gradient-based algorithms
may lead to a local optimum. In this paper, we begin by limiting agents to a
finite set of feasible positions. An advantage of this formulation is that it
assists us in eliminating obviously bad initial conditions for any
gradient-based method. An additional advantage comes from the fact that we can
show our coverage objective function to be monotone submodular, therefore, a
suboptimal solution obtained by the greedy algorithm can achieve a performance
ratio $L(N)\geq1-\frac{1}{e}$, where $N$ is the number of agents in the
system. The idea of exploiting the submodularity of the objective function in
optimization problems has been used in the literature, e.g., in sensor
placement \cite{krause2008near,krause2008efficient} and the maximum coverage
problem mentioned above, whereas a total backward curvature of string
submodular functions is proposed in \cite{zhang2016string} and a total
curvature $c_{k}$ for the $k$-batch greedy algorithm is proposed in
\cite{liu2016performance} in order to derive bounds for related problems.

Our goal in this paper is to derive a tighter lower bound, i.e., to increase
the ratio $L(N)$ by further exploiting the structure of our objective
function. In particular, we make use of the \emph{total curvature}
\cite{conforti1984submodular} and the \emph{elemental curvature}
\cite{wang2016approximation} of the objective function and show that these can
be explicitly derived and lead to new and tighter lower bounds. Moreover, we
show that the tightness of the lower bounds obtained through the total
curvature and the elemental curvature respectively is \emph{complementary}
with respect to the sensing capabilities of the agents. In other words, when
the sensing capabilities are weak, one of the two bounds is tight and when the
sensing capabilities are strong, the other bound is tight. Thus, regardless of
the sensing properties of our agents, we can always determine a lower bound
tighter than $L(N)=1-\frac{1}{e}$ and, in some cases very close to $1$,
implying that the greedy algorithm solution can be guaranteed to be
near-globally optimal.

Another contribution of the paper is to add a final step to the optimal
coverage process, after obtaining the greedy algorithm solution and evaluating
the associated lower bound with respect to the global optimum. Specifically,
we relax the set of allowable agent positions in the mission space from the
imposed discrete set and use the solution of the greedy algorithm as an
initial condition for the distributed gradient-based algorithm in
\cite{Sun2014}. We refer to this as the \emph{Greedy-Gradient Algorithm} (GGA)
which is applicable to the original coverage problem.

The remainder of this paper is organized as follows. The optimal coverage
problem is formulated in Sec.~\ref{II}. In Sec. III, we review key elements of
the submodularity theory and show that how to apply it to the optimal coverage
problem. The GGA for the optimal coverage problem is presented in Sec. IV. In
Sec. V, we provide simulation examples to show how the algorithm works and can
provide significantly better performance compared to earlier results reported
in \cite{Sun2014} .

\section{OPTIMAL COVERAGE PROBLEM FORMULATION}

\label{II} 
\begin{figure}[ptb]
\centering
\includegraphics[
width=0.375\textwidth]
{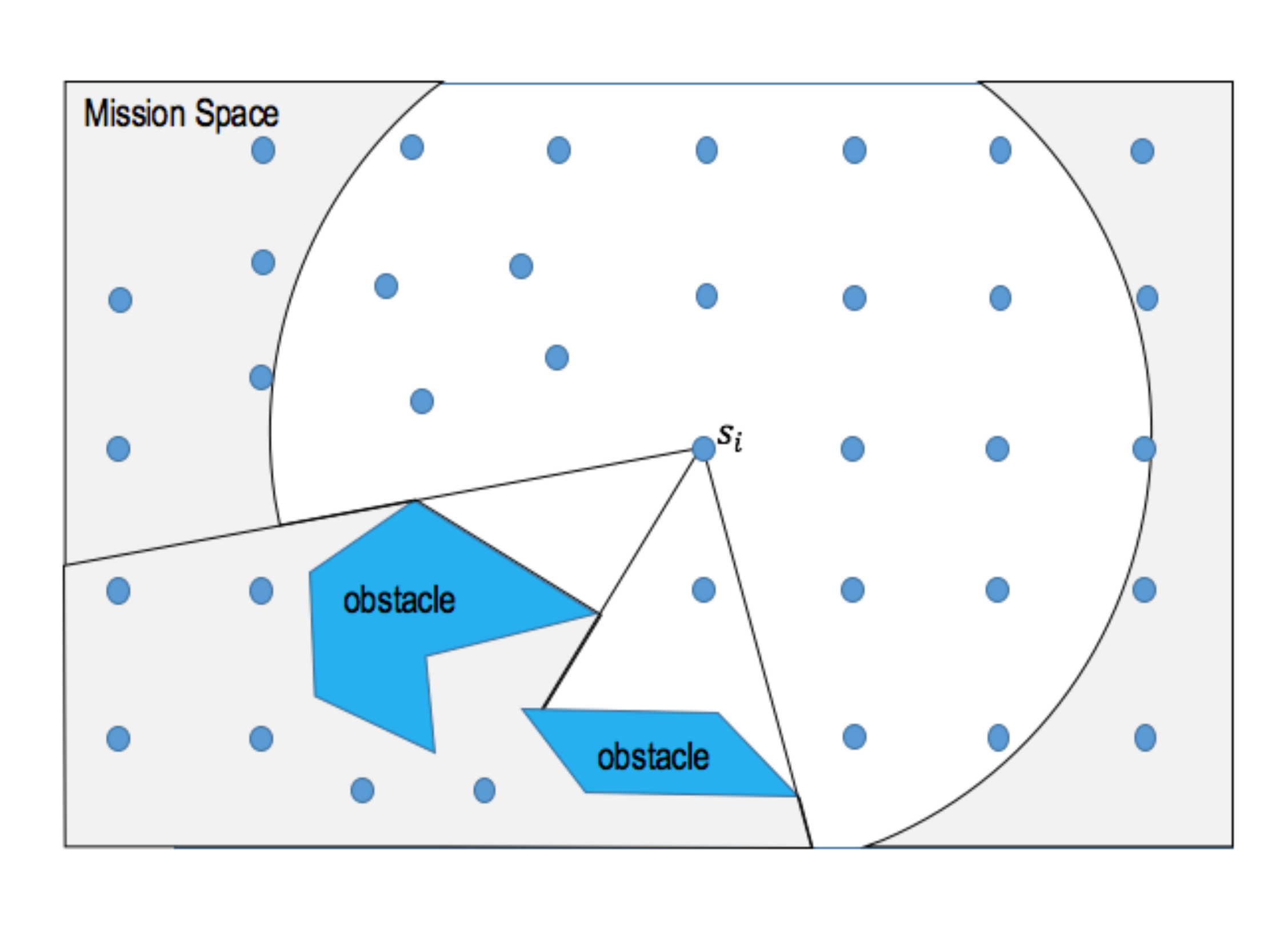} \caption{Mission space example, $F^{D}$ consists of the
blue dots }%
\label{fig:DiscreteFeasibleSpace}%
\end{figure}

We begin by reviewing the basic coverage problem presented in
\cite{caicedo2008performing,CM2004,Minyi2011}. A \textit{mission space}
$\Omega\subset\mathbb{R}^{2}$ is modeled as a non-self-intersecting polygon,
i.e., a polygon such that any two non-consecutive edges do not intersect.
Associated with $\Omega$, we define a function $R(x):\Omega\rightarrow
\mathbb{R}$ to characterize the probability of event occurrences at the
location $x\in\Omega$. It is referred to as \textit{event density} satisfying
$R(x)\geq0$ for all $x\in\Omega$ and $\int_{\Omega}R(x)dx<\infty$. The mission
space may contain obstacles modeled as $m$ non-self-intersecting polygons
denoted by $M_{j}$, $j=1,\ldots,m$, which block the movement as well as the
sensing range of an agent. The interior of $M_{j}$ is denoted by
$\mathring{M_{j}}$ and the overall \textit{feasible space} is $F=\Omega
\setminus(\mathring{M_{1}}\cup\ldots\cup\mathring{M_{m}})$, i.e., the space
$\Omega$ excluding all interior points of the obstacles. There are $N$ agents
in the mission space and their positions are defined by a vector
$\mathbf{s}=(s_{1},\ldots,s_{N})$ with $s_{i}\in F^{D}$, $i=1,\ldots,N$, where
$F^{D}=\{f_{1},\ldots,f_{n}\}$ is a discrete set of feasible positions with
cardinality $n$. We assume that $s_{i}\neq s_{j}$ for any two distinct agents
$i$ and $j$. Figure~\ref{fig:DiscreteFeasibleSpace} shows a mission space with
two obstacles and an agent located at $s_{i}$.

In the coverage problem, agents are sensor nodes. We assume that each node has
a bounded sensing range captured by the \textit{sensing radius} $\delta_{i}$.
Thus, the sensing region of node $i$ is $\Omega_{i}=\{x:d_{i}(x)\leq\delta
_{i}\}$, where $d_{i}(x)=\Vert x-s_{i}\Vert$. The presence of obstacles
inhibits the sensing ability of a node, which motivates the definition of a
\textit{visibility set} $V(s_{i})\subset F$. A point $x\in F$ is
\textit{visible} from $s_{i}\in F$ if the line segment defined by $x$ and
$s_{i}$ is contained in $F$, i.e., $\eta x+(1-\eta)s_{i}\in F$ for all
$\eta\in\lbrack0,1]$, and $x$ is within the sensing range of $s_{i}$, i.e.
$x\in\Omega_{i}$. Then, $V(s_{i})=\Omega_{i}\cap\{x:\eta x+(1-\eta)s_{i}\in
F\,\text{for all }\eta\in\lbrack0,1]\}$ is a set of points in $F$ which are
visible from $s_{i}$. We also define $\bar{V}(s_{i})=F\setminus V(s_{i})$ to
be the \emph{invisibility set} from $s_{i}$, e.g., the grey area in Fig.
\ref{fig:DiscreteFeasibleSpace}. A sensing model for node $i$ is given by the probability that
sensor $i$ detects an event occurrence at $x\in V(s_{i})$, denoted by
$p_{i}(x,s_{i})$. We assume that $p_{i}(x,s_{i})$ can be expressed as a
function of $d_{i}(x)=\Vert x-s_{i}\Vert$ and is monotonically decreasing and
differentiable. An example of such a function is%
\begin{equation}
p_{i}(x,s_{i})=\exp({-\lambda_{i}\Vert x-s_{i}\Vert}), \label{pi}%
\end{equation}
where $\lambda_{i}$ is a \emph{sensing decay} factor. For points that are
invisible to node $i$, the detection probability is zero. Thus, the overall
\textit{sensing detection probability}, denoted by $\hat{p}_{i}(x,s_{i})$, is
defined as%
\begin{equation}
\hat{p}_{i}(x,s_{i})=%
\begin{cases}
p_{i}(x,s_{i}) & \text{if}\quad x\in V(s_{i}),\\
0 & \text{if}\quad x\in\bar{V}(s_{i}),
\end{cases}
\label{phat}%
\end{equation}
which is not a continuous function of $s_{i}$. Note that $V(s_{i}%
)\subset\Omega_{i}=\{x:d_{i}(x)\leq\delta_{i}\}$ is limited by the sensing
range of agents $\delta_{i}$ and that the overall sensing detection
probability of agents is determined by the sensing range $\delta_{i}$ as well
as sensing decay rate $\lambda_{i}$. Then, the \textit{joint detection
probability} that an event at $x\in\Omega$ is detected by the $N$ nodes is
given by%
\begin{equation}
{\label{jointP}}P(x,\mathbf{s})=1-\prod_{i=1}^{N}[1-\hat{p}_{i}(x,s_{i})],
\end{equation}
where we assume that detection probabilities of different sensors are
independent. Assume that $R(x)=0$ for $x\notin F$. The optimal coverage
problem can be expressed as follows:%
\begin{equation}%
\begin{split}
\max_{\mathbf{s}}\text{ }  &  H(\mathbf{s}) =\int_{\Omega}R(x)P(x,\mathbf{s}%
)dx\\
\mathrm{s.t.} \quad &  \mathbf{s}\in\mathcal{I}%
\end{split}
\label{covcontrolH}%
\end{equation}
where $\mathcal{I}=\{S\subseteq F^{D}:|S|\leq N\}$ is a collection of subsets
of $F^{D}$ and $|S|$ denotes the cardinality of set $S$. We emphasize again
that $H(\mathbf{s})$ is not convex (concave) even in the simplest possible
problem setting.

\section{SUBMODULARITY THEORY APPLIED TO THE OPTIMAL COVERAGE PROBLEM}

A naive method to find the global optimum of (\ref{covcontrolH}) is the
brute-force search. The time complexity is $n!/(N!(n-N)!)$ by choosing $N$
agent positions from $n$ feasible positions. The brute-force method may not
generate quality solutions in a reasonable amount of time when $n$ and $N$ are
large. In this section, we will introduce the basic elements of submodularity
theory and apply it to the optimal coverage problem. We will show that our
objective function $H(\mathbf{s})$ in (\ref{covcontrolH}) is \emph{monotone
submodular}, therefore, we can apply basic results from submodularity theory
which hold for this class of functions. According to this theory, the greedy
algorithm (described in Section I and shown in \textbf{Algorithm 1}) produces
a guaranteed performance in polynomial time. The time complexity of the greedy
algorithm is $O(nN)$. When $n$ is given, it is $O(N)$, which is linear in the
number of agents.

\subsection{Monotone Submodular Coverage Metric}

A submodular function is a set function whose value has the diminishing
returns property. The formal definition of submodularity is given as follows.

\textbf{Definition 1: } Given a ground set $Y=\{y_{1},\ldots,y_{n}$\} and its
power set $2^{Y}$, a function $f:2^{Y}\rightarrow\mathbb{R}$ is called
\textit{submodular} if for any $S,T\subseteq Y$,%
\begin{equation}
f(S\cup T)+f(S\cap T)\leq f(S)+f(T).\label{def1}%
\end{equation}
If, additionally, $f(S)\leq f(T)$ whenever $S\subseteq T$, we say that $f$ is
\textit{monotone submodular}. An equivalent definition, which better reflects
the diminishing returns property, is given below, where the proof of
equivalence can be found in Appendix~\ref{app1}.

\textbf{Definition 2: } For any sets $S,T\subseteq Y$ with $S\subseteq T$ and
any $y\in Y\setminus T$, we have%
\begin{equation}
f(S\cup\{y\})-f(S)\geq f(T\cup\{y\})-f(T). \label{def2}%
\end{equation}
Intuitively, the incremental increase of the function is larger when an
element is added to a small set than to a larger set. In what follows, we will
use the second definition.

A general form of the submodular maximization problem is
\begin{equation}%
\begin{split}
\max\quad &  f(S)\\
\mathrm{s.t.}\quad S  &  \in\mathcal{I}%
\end{split}
\label{generalSubmodularFun}%
\end{equation}
where $\mathcal{I}$ is a non-empty collection of subsets of a finite set $Y$.
$\mathcal{M}=(Y,\mathcal{I}),\mathcal{I}\subseteq2^{Y}$ is
\textit{independent} if, for all $B\in\mathcal{I}$, any set $A\subseteq B$ is
also in $\mathcal{I}$. Furthermore, if for all $A\in\mathcal{I},$
$B\in\mathcal{I}$, $|A|<|B|$, there exists a $j\in B\setminus A$ such that
$A\cup\{j\}\in\mathcal{I}$, then $\mathcal{M}$ is called a \textit{matroid}.
Moreover, $\mathcal{M}=(Y,\mathcal{I})$ is called \textit{uniform matroid} if
$\mathcal{I}=\{S\subseteq Y:|S|\leq N\}$.

The following theorem establishes the fact that the objective function
$H(\mathbf{s})$ in (\ref{covcontrolH}) is monotone submodular, regardless of
the obstacles that may be present in the mission space. This will allow us to
apply results that quantify a solution obtained through the greedy algorithm
relative to the global optimum in (\ref{covcontrolH}).

\textbf{Theorem 1: } $H (\mathbf{s})$ is monotone submodular, i.e.,%
\[
H(S \cup\{s_{k}\}) - H(S) \geq H(T \cup\{s_{k}\}) - H(T)
\]
and%
\[
H(S) \leq H(T)
\]
for any $S,T \subseteq F^{D}$ with $S \subseteq T$ and $s_{k} \in F^{D}
\setminus T $.

\proof Let $S$ and $T$, such that $S\subseteq T\subseteq F^{D}$, be two agent
position vectors. Since $S\subseteq T$ and $0\leq1-\hat{p}_{i}(x,s_{i})\leq1$
for any $s_{i}\in F^{D}$, we have
\begin{equation}
\prod_{s_{i}\in S}[1-\hat{p}_{i}(x,s_{i})]\geq\prod_{s_{i}\in T}[1-\hat{p}%
_{i}(x,s_{i})] \label{subinq}%
\end{equation}
for all $x\in\Omega$. In addition, $H(S\cup\{s_{k}\})$ can be written as%
\[%
\begin{split}
&  H(S\cup\{s_{k}\})\\
=  &  \int_{\Omega}R(x)\left\{  1-[1-\hat{p}_{k}(x,s_{k})]\prod_{s_{i}\in
S}[1-\hat{p}_{i}(x,s_{i})]\right\}  dx\\
=  &  \int_{\Omega}R(x)\left\{  1-\prod_{s_{i}\in S}[1-\hat{p}_{i}%
(x,s_{i})]\right\}  dx\\
&  +\int_{\Omega}R(x)\hat{p}_{k}(x,s_{k})\prod_{s_{i}\in S}[1-\hat{p}%
_{i}(x,s_{i})]dx.
\end{split}
\]
The difference between $H(S)$ and $H(S\cup\{s_{k}\})$ is given by%
\begin{equation}%
\begin{split}
&  H(S\cup\{s_{k}\})-H(S)\\
&  =\int_{\Omega}R(x)\hat{p}_{k}(x,s_{k})\prod_{s_{i}\in S}[1-\hat{p}%
_{i}(x,s_{i})]dx.
\end{split}
\label{differenceOneMore}%
\end{equation}
Using the same derivation for $T$, we can obtain
\begin{equation}%
\begin{split}
&  H(T\cup\{s_{k}\})-H(T)\\
&  =\int_{\Omega}R(x)\hat{p}_{k}(x,s_{k})\prod_{s_{i}\in T}[1-\hat{p}%
_{i}(x,s_{i})]dx.
\end{split}
\label{differenceOneMoreT}%
\end{equation}
From (\ref{differenceOneMore}) and (\ref{differenceOneMoreT}), the difference
between $H(S\cup\{s_{k}\})-H(S)$ and $H(T\cup\{s_{k}\})-H(T)$ is%
\begin{equation}%
\begin{split}
&  \left[  H(S\cup\{s_{k}\})-H(S)\right]  -\left[  H(T\cup\{s_{k}%
\})-H(T)\right] \\
=  &  \int_{\Omega}R(x)\hat{p}_{k}(x,s_{k})\prod_{s_{i}\in S}[1-\hat{p}%
_{i}(x,s_{i})]dx\\
&  -\int_{\Omega}R(x)\hat{p}_{k}(x,s_{k})\prod_{s_{i}\in T}[1-\hat{p}%
_{i}(x,s_{i})]dx.
\end{split}
\end{equation}

Using (\ref{subinq}), it follows that the difference $[H(S\cup\{s_{k}%
\})-H(S)]-[H(T\cup\{s_{k}\})-H(T)]\geq0$. Therefore, from Definition 2,
$H(\mathbf{s})$ is submodular.

Next, we prove that $H(\mathbf{s})$ is monotone, i.e., $H(S)\leq H(T)$.
Subtracting $H(T)$ from $H(S)$ yields%
\[%
\begin{split}
&  H(S)-H(T)\\
= &  \int_{\Omega}R(x)\left\{  1-\prod_{s_{i}\in S}[1-\hat{p}_{i}%
(x,s_{i})]\right\}  dx\\
&  -\int_{\Omega}R(x)\left\{  1-\prod_{s_{i}\in T}[1-\hat{p}_{i}%
(x,s_{i})]\right\}  dx\\
= &  \int_{\Omega}R(x)\left\{  \prod_{s_{i}\in T}[1-\hat{p}_{i}(x,s_{i}%
)]-\prod_{s_{i}\in S}[1-\hat{p}_{i}(x,s_{i})]\right\}  dx.
\end{split}
\]
Using (\ref{subinq}), we have $H(S)-H(T)\leq0$. Therefore, $H(\mathbf{s})$ is
a monotone submodular function. \hfill$\blacksquare$

\subsection{Greedy Algorithm and Lower Bounds}

Finding the optimal solution to (\ref{generalSubmodularFun}) is in general
NP-hard. The following greedy algorithm is usually used to obtain a feasible
solution for (\ref{generalSubmodularFun}). The basic idea of the greedy
algorithm is to add an agent which can maximize the value of the objective
function at each iteration.
\renewcommand{\algorithmicrequire}{\textbf{Input:}}
\renewcommand{\algorithmicensure}{\textbf{Output:}} \begin{algorithm}[H]
\caption{Greedy Algorithm}
\label{alg: greedyAlgorithm}
\begin{algorithmic}[1]
\REQUIRE Submodular function $f(S)$\\
Cardinality constraint $N$
\ENSURE Set $S$\\
\hspace{-0.6cm}\textbf{Initialization:} $S \leftarrow \emptyset$, $i \leftarrow 0$
\WHILE{$i \leq N$}
\STATE $s_i^{\ast} = \argmax_{s_i\in Y\setminus S} {f(S \cup \{s_i\} )}$
\STATE $S \leftarrow S \cup\{s_i^{\ast}\}$
\STATE $i \leftarrow i+1 $
\ENDWHILE
\RETURN $S$
\end{algorithmic}
\end{algorithm}

In the following analysis, we assume that $f$ is a monotone submodular
function satisfying $f(\emptyset)=0$ and $\mathcal{M}=(Y,\mathcal{I})$ is a
uniform matroid. We will use the definition
\[
L(N)=\frac{f}{f^{\ast}}%
\]
from Section I, where $f^{\star}$ is the global optimum of
(\ref{generalSubmodularFun}) and $f$ is a feasible solution obtained by
Algorithm~\ref{alg: greedyAlgorithm}. Then, as shown in
\cite{nemhauser1978analysis}, a lower bound of $L(N)$ is $1-1/e$ .

Next, we consider the \textit{total curvature}%
\begin{equation}
c=\max_{j\in Y}\left[  1-\frac{f(Y)-f(Y\setminus j)}{f(\{j\})}\right]
\label{totalCurvatureDef}%
\end{equation}
introduced in \cite{conforti1984submodular}. Using $c$, the lower bound of
$L(N)$ above is improved to be $T(c,N)$:
\begin{equation}
T(c,N)=\frac{1}{c}\left[  1-(\frac{N-c}{N})^{N}\right]  .
\label{totalCurvLowerBound}%
\end{equation}
where $c\in\lbrack0,1]$, and
\[
T(c,N)\geq1-\frac{1}{e}%
\]
for any $N\geq1$. If $c=1$, the result is the same as the bound obtained in
\cite{nemhauser1978analysis}, \cite{fisher1978analysis}.

In addition, we consider the \textit{elemental curvature}%
\begin{equation}
\alpha=\max_{S\subset Y,i,j\in Y\setminus S,i\neq j}\frac{f(S\cup
\{i,j\})-f(S\cup\{j\})}{f(S\cup\{i\})-f(S)}, \label{elemenCurvatureDef}%
\end{equation}
based on which the following bound is obtained:%
\begin{equation}
E(\alpha,N)=1-\left(  \frac{\alpha+\ldots+\alpha^{N-1}}{1+\alpha+\ldots
+\alpha^{N-1}}\right)  ^{N} \label{eleCurvLowerBound}%
\end{equation}
and it is shown in \cite{wang2016approximation} that $L(N)\geq E(\alpha,N)$.
Note that $E(\alpha,N)$ can be simplified as follows:%
\begin{equation}
E(\alpha,N)=%
\begin{cases}
1-(\frac{N-1}{N})^{N}, & \text{when }\alpha=1;\\
1-(\frac{\alpha-\alpha^{N}}{1-\alpha^{N}})^{N}, & \text{when }0\leq\alpha<1.
\end{cases}
\end{equation}

If both bounds $T(c,N)$ and $E(\alpha,N)$ can be calculated, then the larger
one will be the lower bound $L(N)$, defined as
\begin{equation}
L(N)=\max\{T(c,N),E(\alpha,N)\}.\label{Lmax}%
\end{equation}
Accordingly, we have $f(S)\geq L(N)f(S^{\ast})$, where $S^{\star}$ is the
global optimum set, and $S$ is the set obtained by
Algorithm~\ref{alg: greedyAlgorithm}.

\begin{figure}[ptb]
\centering
\includegraphics[
height=2.8in,
width=3in]
{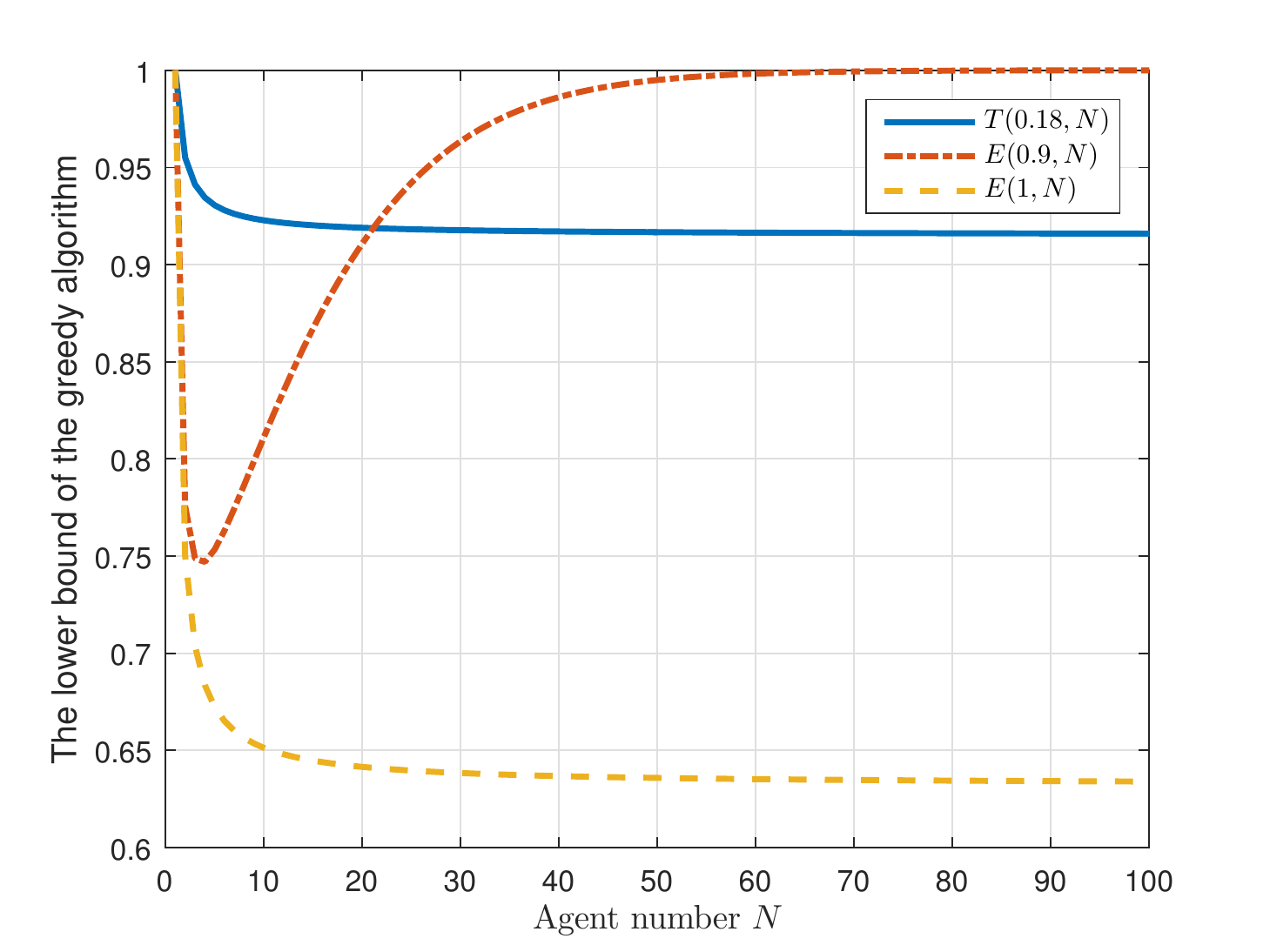} \caption{ $T(c,N)$ and $E(\alpha,N)$ as a function of the
number of agents $N$}%
\label{fig:CBoundNvary}%
\end{figure}

Figure \ref{fig:CBoundNvary} shows the dependence of $T(c,N)$ and
$E(\alpha,N)$ on the number of agents $N$ for some specific values of $c$ and
$\alpha$ (as shown in the figure). Clearly, if $c<1$ and $\alpha<1$, then
$L(N)$ in (\ref{Lmax}) is much tighter than $1-\frac{1}{e}$.

\subsection{Curvature Information Calculation}

In this subsection, we will derive the concrete form of the total curvature
$c$ and the elemental curvature $\alpha$ in the context of coverage problems.
For notational convenience, $\hat{p}_{i}(x,s_{i})$ is used without its
arguments as long as this dependence is clear from the context.

Recall that $F^{D}$ is the set of feasible agent positions. We can obtain from
(\ref{covcontrolH}):
\[
{\label{objFullSep}}%
\begin{split}
H(F^{D}) &  =\int_{\Omega}R(x)\left[  1-\prod_{i=1}^{n}(1-\hat{p}_{i})\right]
dx\\
&  =\int_{\Omega}R(x)\left[  1-(1-\hat{p}_{j})\prod_{i=1,i\neq j}^{n}%
(1-\hat{p}_{i})\right]  dx,
\end{split}
\]
and%
\[
{\label{objFullSepMinus}}H(F^{D}\setminus\{s_{j}\})=\int_{\Omega}R(x)\left[
1-\prod_{i=1,i\neq j}^{n}(1-\hat{p}_{i})\right]  dx.
\]
The difference between $H(F^{D})$ and $H(F^{D}\setminus\{s_{j}\})$ is%
\begin{equation}
H(F^{D})-H(F^{D}\setminus\{s_{j}\})=\int_{\Omega}R(x)\hat{p}_{j}%
\prod_{i=1,i\neq j}^{n}[1-\hat{p}_{i}]dx.\label{difference}%
\end{equation}
When there is only one agent $s_{j}$, the objective function is
\begin{equation}
H(s_{j})=\int_{\Omega}R(x)\hat{p}_{j}dx.\label{OneAgent}%
\end{equation}
Combining (\ref{totalCurvatureDef}), (\ref{difference}) and (\ref{OneAgent}),
we obtain
\begin{equation}
c=\max_{s_{j}\in F^{D}}\left[  1-\frac{\int_{\Omega}R(x)\hat{p}_{j}%
\prod_{i=1,i\neq j}^{n}[1-\hat{p}_{i}]dx}{\int_{\Omega}R(x)\hat{p}_{j}%
dx}\right]  .
\end{equation}
\textbf{Remark 1} If the sensing capabilities of agents are weak, that is,
$\hat{p}_{i}$ is small for most parts in the mission space, then
$\prod_{i=1,i\neq j}^{n}(1-\hat{p}_{i})$ is, in turn, close to 1, which leads
to a small value of $c$. It follows from (\ref{totalCurvLowerBound}) that the
lower bound $T(c,N)$ is a monotonically decreasing function of $c$ and
approaches 1 near $c=0$. This implies that the solution of the greedy
algorithm is very close to the global optimum when the sensing capabilities
are weak.

Next, we calculate the elemental curvature $\alpha$. From
(\ref{differenceOneMore}), the difference between $H(S)$ and $H(S\cup
\{s_{k}\})$ is%
\begin{equation}
H(S\cup\{s_{k}\})-H(S)=\int_{\Omega}R(x)\hat{p}_{k}(x)\prod_{s_{i}\in
S}[1-\hat{p}_{i}]dx.
\end{equation}
Using the same derivation, we can obtain%
\begin{equation}%
\begin{split}
&  H(S\cup\{s_{j},s_{k}\})-H(S\cup\{s_{j}\})\\
&  \qquad=\int_{\Omega}R(x)\hat{p}_{k}(1-\hat{p}_{j})\prod_{s_{i}\in S}%
[1-\hat{p}_{i}]dx.
\end{split}
\end{equation}
The elemental curvature in (\ref{elemenCurvatureDef}) can then be calculated
by
\begin{equation}%
\begin{split}
\alpha &  =\max_{S,s_{j},s_{k}}\frac{H(S\cup\{s_{j},s_{k}\})-H(S\cup
\{s_{j}\})}{H(S\cup\{s_{k}\})-H(S)}\\
&  =\max_{S,s_{j},s_{k}}\frac{\int_{\Omega}R(x)\hat{p}_{k}(1-\hat{p}_{j}%
)\prod_{s_{i}\in S}[1-\hat{p}_{i}]dx}{\int_{\Omega}R(x)\hat{p}_{k}\prod
_{s_{i}\in S}[1-\hat{p}_{i}]dx}\\
&  =\max_{S,s_{j},s_{k}}1-\frac{\int_{\Omega}R(x)\hat{p}_{k}\hat{p}_{j}%
\prod_{s_{i}\in S}[1-\hat{p}_{i}]dx}{\int_{\Omega}R(x)\hat{p}_{k}\prod
_{s_{i}\in S}[1-\hat{p}_{i}]dx}\\
&  =1-\min_{s_{j},x\in\Omega}\hat{p}_{j}(x,s_{j}).
\end{split}
\label{elemenCurvatureCoverage}%
\end{equation}
\textbf{Remark 2} Observe that the elemental curvature turns out to be
determined by a single agent. If there exists a pair $(x,s_{j})$ such that
$x\in\bar{V}(s_{j})$ in (\ref{phat}), then $\hat{p}_{j}(x,s_{j})=0$ and
$\alpha=1$. This may happen when there are obstacles in the mission space or
the sensing capabilities of agents are weak (e.g., the sensing range is small
or the sensing decay rate is large). On the other hand, if the sensing
capabilities are so strong that $\hat{p}_{j}(x,s_{j})\neq0$ for any $x\in
F,s_{j}\in F^{D}$, then $\alpha<1$. In addition, $E(\alpha,N)$ is a
monotonically decreasing function of $\alpha$.

\begin{figure}[ptb]
\centering
\includegraphics[
height=2.8in,
width=3in]
{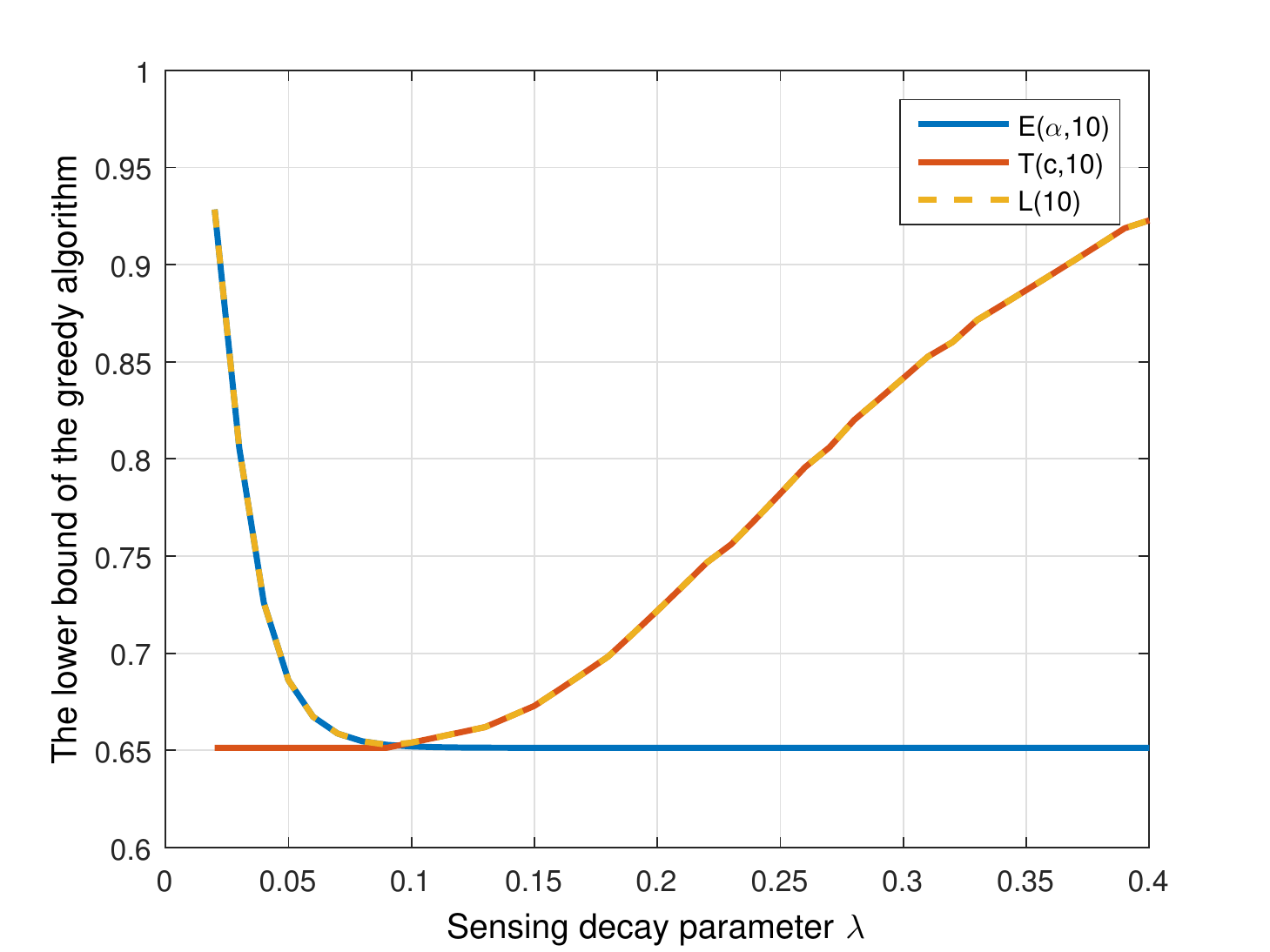} \caption{Lower bound $L(10)$ as a function of the sensing
decay rate of agents}%
\label{fig:Bound}%
\end{figure}

\begin{figure}
\centering
\includegraphics[
height=2.8in,
width=3in]
{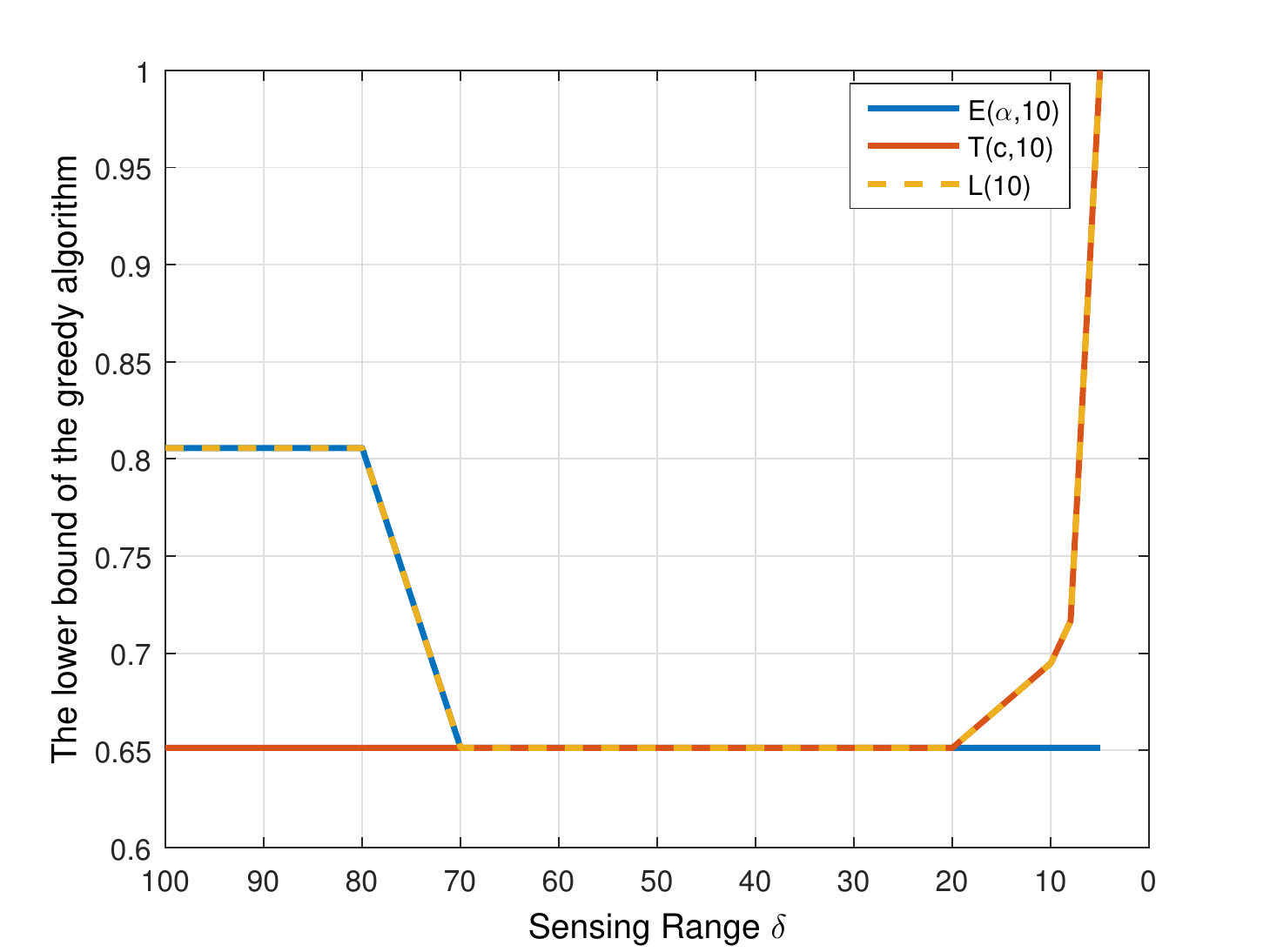} \caption{Lower bound $L(10)$ as a
function of the sensing range of agents}%
\label{fig:Bound2}%
\end{figure}

An interesting conclusion from this analysis is that $T(c,N)$ and
$E(\alpha,N)$ are \emph{complementary} with respect to the sensing
capabilities of sensors. From Remark 1, $T(c,N)$ is large when the sensing
capabilities are weak, while from Remark 2, $E(\alpha,N)$ is large when the
sensing capabilities are strong. This conclusion is graphically depicted in
Figs. \ref{fig:Bound} and \ref{fig:Bound2} (where sensing capability varies
from strong to weak). In Fig. \ref{fig:Bound}, $E(\alpha,N)$ and $T(c,N)$ have
been evaluated for $N=10$ and $\delta=80$ as a function of one of the measures
of sensing capability, the sensing decay rate $\lambda$ in (\ref{pi}),
assuming all agents have the same sensing capabilities. One can see that for
small values of $\lambda$, the bound $E(\alpha,10)$ is close to 1 and
dominates both $T(c,10)$ and the well-known bound $1-\frac{1}{e}$. Beyond a
critical value of $\lambda$, it is $T(c,10)$ that dominates and approaches 1
for large values of $\lambda$. Figure \ref{fig:Bound2} shows a similar
behavior when $T(c,N)$ and $E(\alpha,N)$ are evaluated for $N=10$ and $\lambda=0.03$ as a function of the
other measure of sensing capability, the sensing range $\delta$. When the
sensing range exceeds the distance of the diagonal of the mission space, there
is no value in further increasing the sensing range and $E(\cdot)$ becomes constant. 
When $\delta>20$, the sensing capabilities are strong and $T(\cdot)$ becomes constant. 
Therefore, both $E(\cdot)$ and $T(\cdot)$ become constant when $\delta$ exceeds corresponding thresholds.
On the other hand, when the sensing range is smaller than some threshold, then
$\alpha=1$, and $E(1,10)=0.6513$. 

Figures \ref{fig:Bound} and \ref{fig:Bound2} also illustrate the trade-off
between the sensing capabilities and the coverage performance guarantee.
Agents with strong capabilities obviously achieve better coverage performance.
On the other hand, one can get a better guaranteed performance as the agents'
capabilities get weaker. Therefore, if one is limited to agents with weak
sensing capabilities in a particular setting, the use of $T(c,N)$ is
appropriate and this trade-off may be exploited.

\section{GREEDY-GRADIENT ALGORITHM}

Thus far, we have restricted agent positions to be selected from the finite
feasible set $f^{D}=\{f_{1},...,f_{n}\}$. In this section, agents are allowed
to be deployed at any feasible point, and the optimal coverage problem
becomes
\begin{equation}%
\begin{split}
\max_{\mathbf{s}}\text{ }  &  H(\mathbf{s})=\int_{\Omega}R(x)P(x,\mathbf{s}%
)dx\\
\mathrm{s.t.}\quad &  s_{i}\in F,i=1,...,N.
\end{split}
\label{covContinueH}%
\end{equation}
We propose a \emph{Greedy-Gradient Algorithm} (GGA) shown in Algorithm
\ref{alg: greedyGradientAlgorithm} to solve this problem. The basic idea is to
use existing gradient-based algorithms with an initial deployment given by the
greedy algorithm (Algorithm 1) to seek better performance. In particular, we
use the distributed gradient-based algorithm developed in \cite{Minyi2011}:
\begin{equation}
s_{i}^{k+1}=s_{i}^{k}+\zeta_{k}\frac{\partial H(\mathbf{s})}{\partial
s_{i}^{k}},\text{ \ }k=0,1,\ldots\label{gradientbasedalgo}%
\end{equation}
where the step size sequence $\{{\zeta_{k}\}}$ is appropriately selected to
ensure convergence of the resulting trajectories for all agents
\cite{Bertsekas95}. The detailed calculation of $\frac{\partial H(\mathbf{s}%
)}{\partial s_{i}^{k}}$ can be found in \cite{Sun2014}.
\renewcommand{\algorithmicrequire}{\textbf{Input:}}
\renewcommand{\algorithmicensure}{\textbf{Output:}} \begin{algorithm}
\caption{Greedy-Gradient Algorithm}
\label{alg: greedyGradientAlgorithm}
\begin{algorithmic}[1]
\REQUIRE  Objective function $H(\textbf{s})$\\
\ENSURE Agent positions $\textbf{s}$\\
\hspace{-0.6cm}\textbf{Initialization:} $\textbf{s}$ given by Greedy Algorithm~\ref{alg: greedyAlgorithm} \\
\WHILE{the stopping criterion is not satisfied}
\STATE Choose a step size $\zeta>0$
\FOR{$i=1,\ldots,N$}
\STATE Determine a searching direction $\frac{\partial H(\mathbf{s})}{\partial s_{i}}$
\STATE Update: $s_{i}\leftarrow s_{i}+\zeta \frac{\partial H(\mathbf{s})}{\partial s_{i}} $
\ENDFOR
\ENDWHILE
\RETURN $\textbf{s}$
\end{algorithmic}
\end{algorithm}The stopping criterion is of the form $\Vert\frac{\partial
H(\mathbf{s})}{\partial s_{i}} \Vert\leq\eta$, where $\eta$ is a small
positive scalar.

\section{SIMULATION RESULTS}

In this section, we illustrate through simulation our analysis and the use of
the GGA (Algorithm 2) for coverage problems in a variety of mission spaces
with and without obstacles. The mission space is a $60\times50$ rectangular
area and the event density function $R(x)$ is assumed to be uniformly
distributed, i.e., we set $R(x)=1$ in (\ref{covcontrolH}). We first compute
the theoretical lower bound $L(N)$ for the case of no obstacles in the mission
space and the number of agents is $N=10$.

\begin{figure}[ptb]
\centering
\begin{subfigure}[t]{1in}
			\centering
			\includegraphics[width=\textwidth]{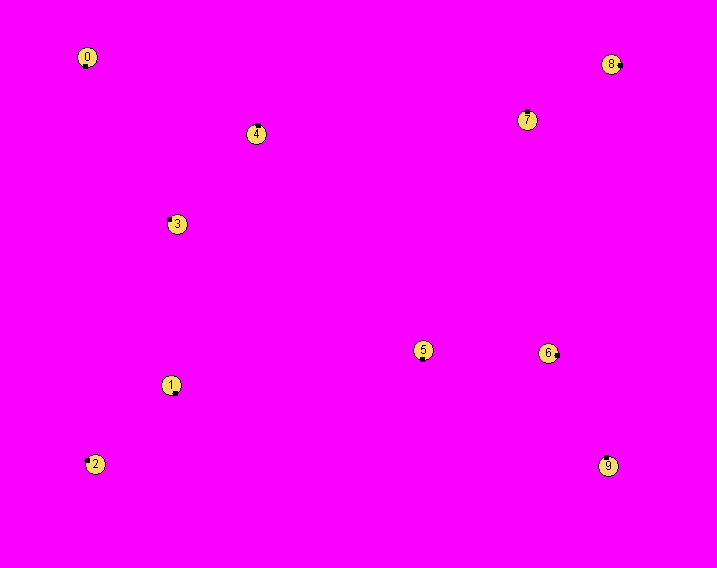}	
			\caption{$H(\textbf{s})=2999.7$}
			\label{fig:GradientP1}
	 \end{subfigure}
\quad\begin{subfigure}[t]{1in}
			\centering
			\includegraphics[width=\textwidth]{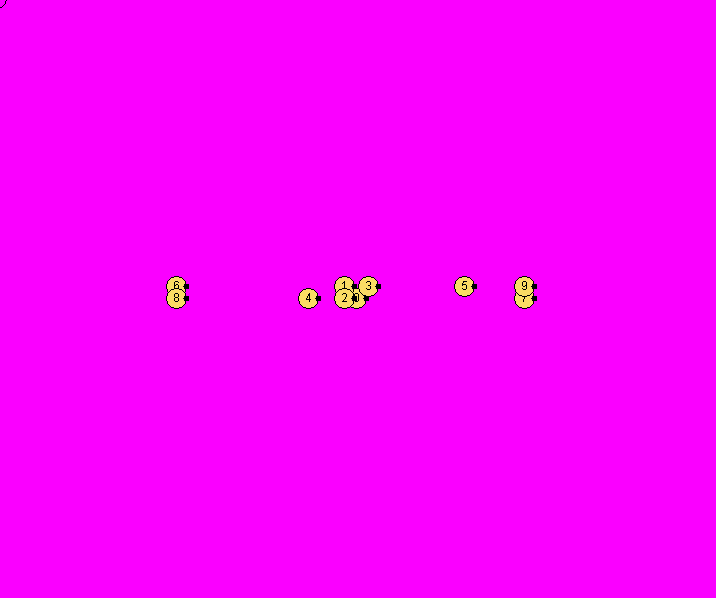}	
			\caption{$H(\textbf{s})=2999.6$}
			\label{fig:GreedyP1}	
	\end{subfigure}
\quad\begin{subfigure}[t]{1in}
			\centering
			\includegraphics[width=\textwidth]{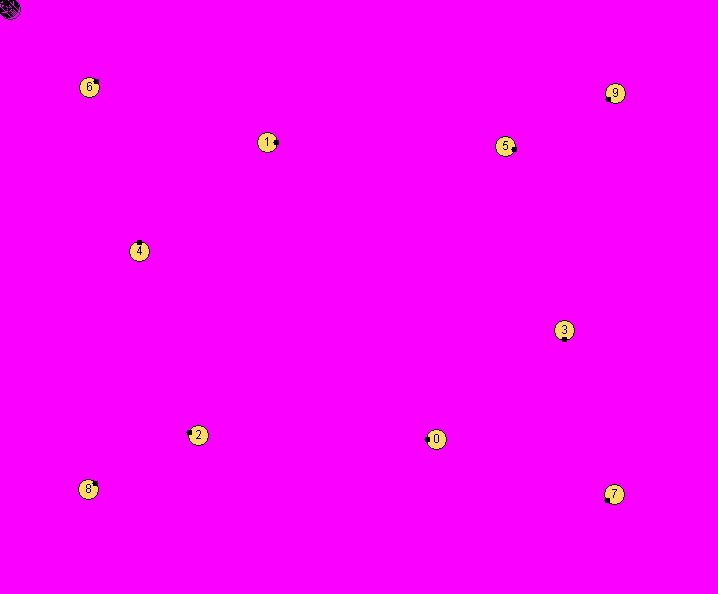}	
			\caption{$H(\textbf{s})=2999.7$}
			\label{fig:GreedyGradientP1}
	 \end{subfigure}
\caption{The decay factor $\lambda=0.02$, and no obstacles in the mission
space}%
\label{fig:NoObsP1}%
\end{figure}

\begin{figure}[ptb]
\centering
\begin{subfigure}[t]{1in}
			\centering
			\includegraphics[width=\textwidth]{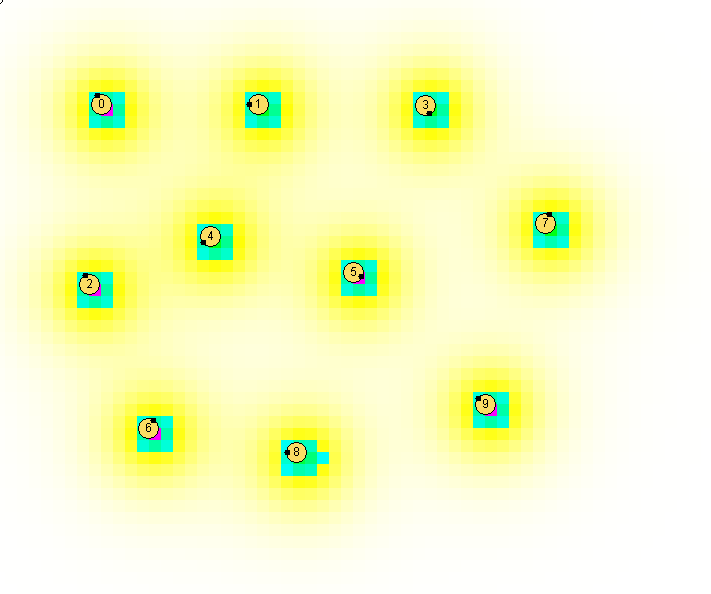}	
			\caption{$H(\textbf{s})=78.3$}
			\label{GradientP3}
	 \end{subfigure}
\quad\begin{subfigure}[t]{1in}
			\centering
			\includegraphics[width=\textwidth]{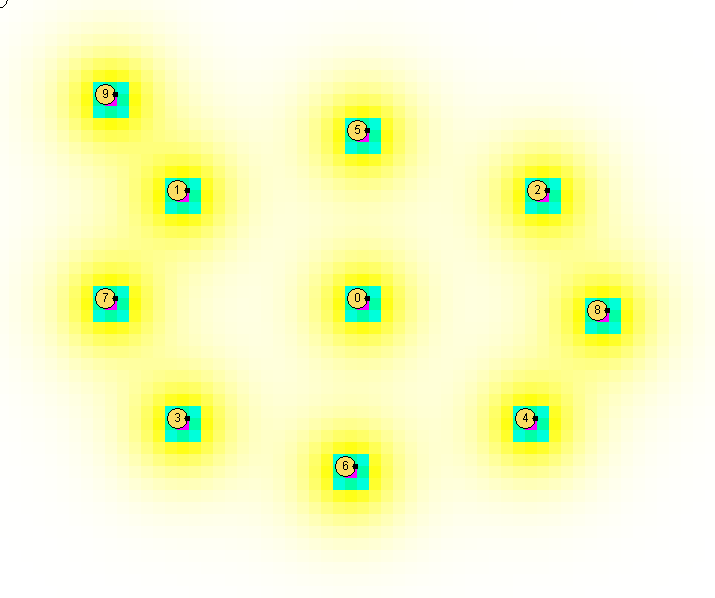}
			\caption{$H(\textbf{s})=78.3$} 	
			\label{fig:GreedyP3}
	\end{subfigure}
\quad\begin{subfigure}[t]{1in}
			\centering
			\includegraphics[width=\textwidth]{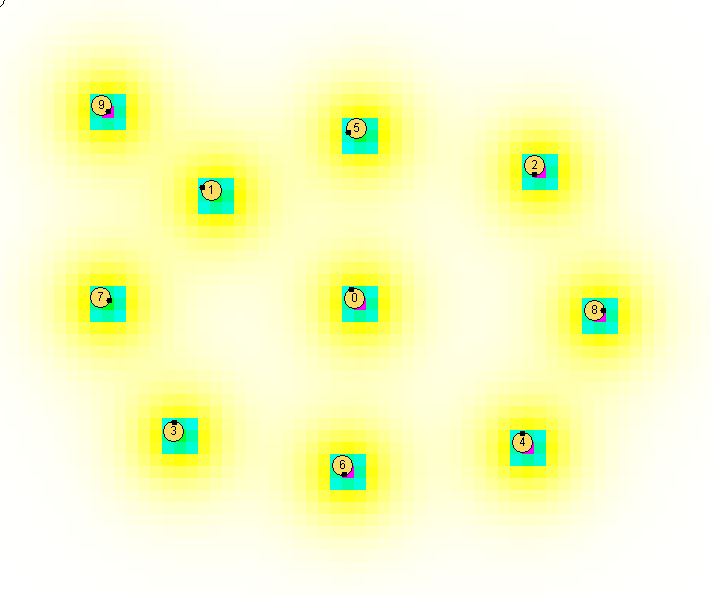}	
			\caption{$H(\textbf{s})=78.3$}
			\label{fig:GreedyGradientP3}
	 \end{subfigure}
\caption{The decay factor $\lambda=0.4$, and no obstacles in the mission
space}%
\label{fig:NoObsP3}%
\end{figure}

\begin{figure}[ptb]
\centering
\begin{subfigure}[t]{1in}
			\centering
			\includegraphics[width=\textwidth]{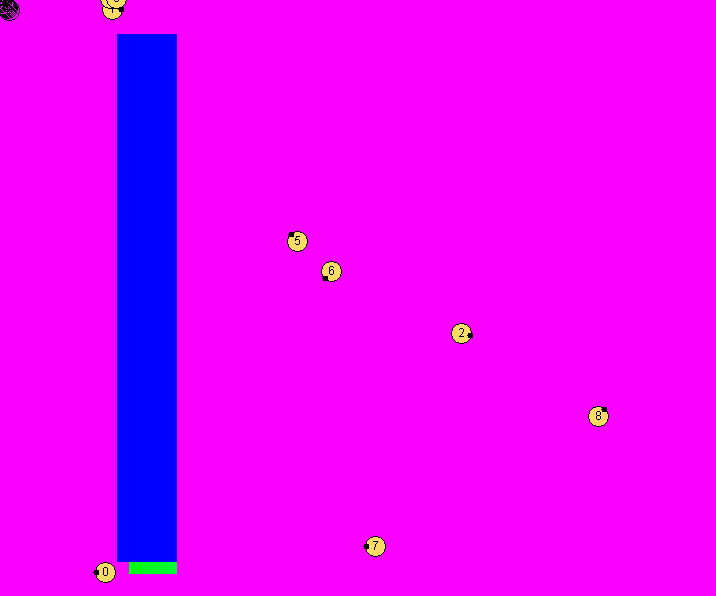}	
			\caption{$H(\textbf{s})=2771.3$}
	 \end{subfigure}
\quad\begin{subfigure}[t]{1in}
			\centering
			\includegraphics[width=\textwidth]{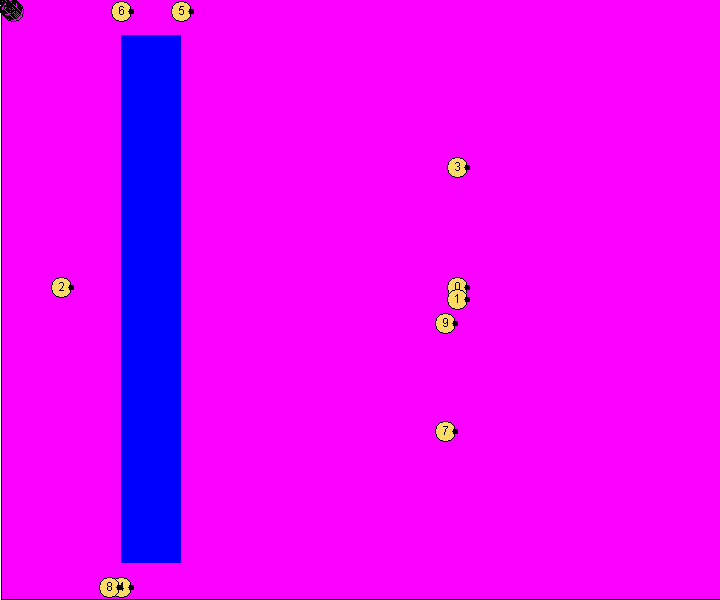}	
			\caption{$H(\textbf{s})=2773.9$} 	
	\end{subfigure}
\quad\begin{subfigure}[t]{1in}
			\centering
			\includegraphics[width=\textwidth]{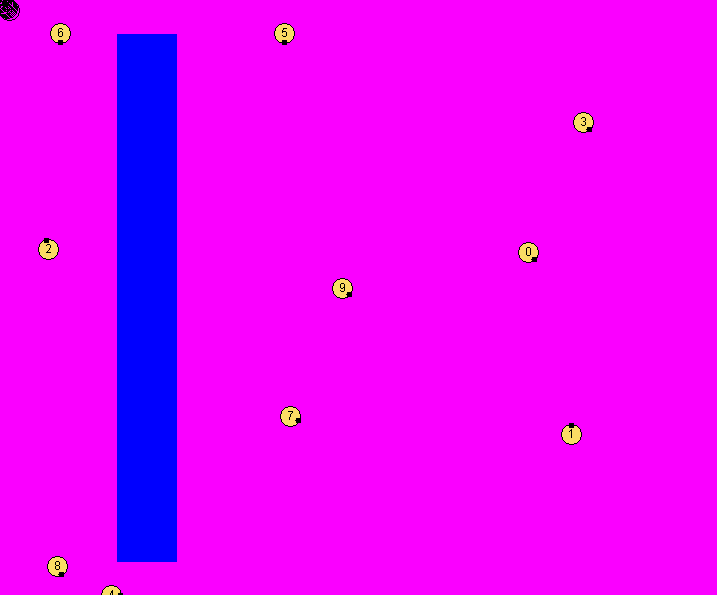}	
			\caption{$H(\textbf{s})=2774.6$}
	 \end{subfigure}
\caption{The decay factor $\lambda=0.02$, and a wall-like obstacle in the
mission space}%
\label{fig:NarrowP1}%
\end{figure}

\begin{figure}[ptb]
\centering
\begin{subfigure}[t]{1in}
			\centering
			\includegraphics[width=\textwidth]{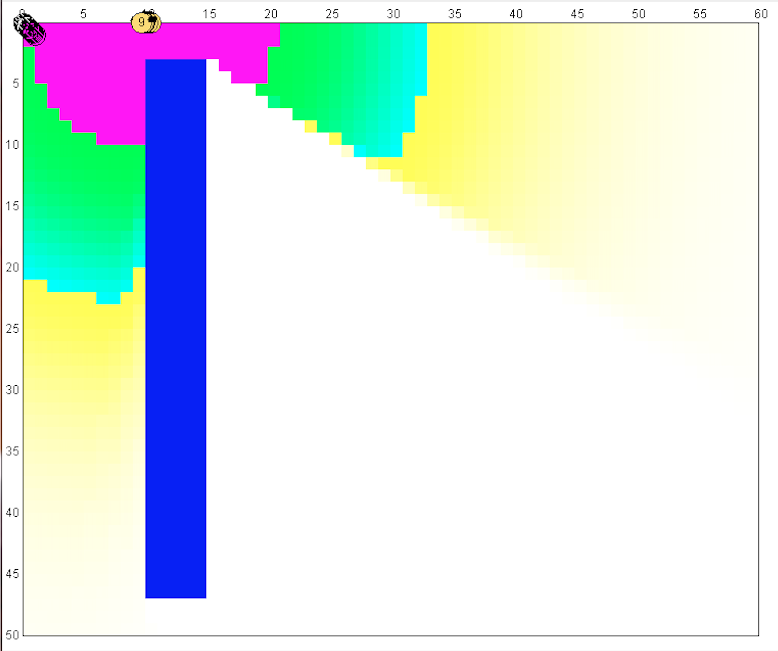}	
			\caption{$H(\textbf{s})=437.1$}
	 \end{subfigure}
\quad\begin{subfigure}[t]{1in}
			\centering
			\includegraphics[width=\textwidth]{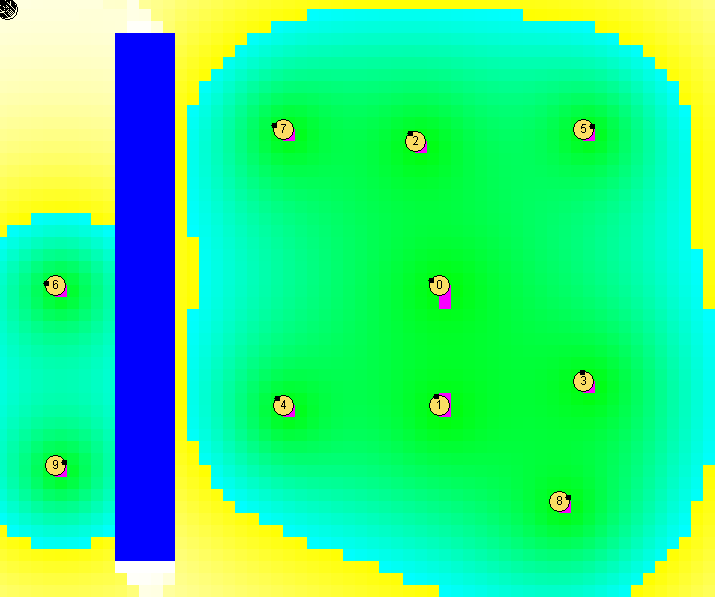}
			\caption{$H(\textbf{s})=1813.3$} 	
	\end{subfigure}
\quad\begin{subfigure}[t]{1in}
			\centering
			\includegraphics[width=\textwidth]{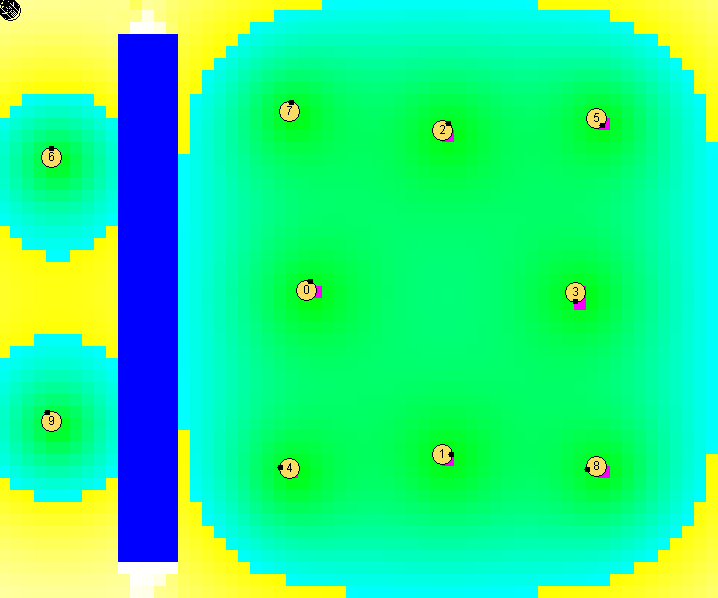}	
			\caption{$H(\textbf{s})=1846.3$}
	 \end{subfigure}
\caption{The decay factor $\lambda=0.12$, and a wall-like obstacle in the
mission space}%
\label{fig:NarrowP2}%
\end{figure}

\begin{figure}[ptb]
\centering
\begin{subfigure}[t]{1in}
			\centering
			\includegraphics[width=\textwidth]{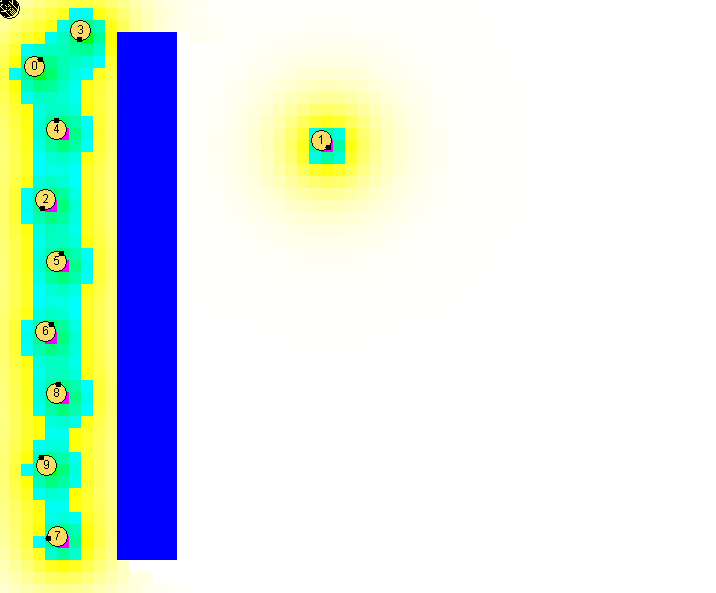}	
			\caption{$H(\textbf{s})=269.6$}
	 \end{subfigure}
\quad\begin{subfigure}[t]{1in}
			\centering
			\includegraphics[width=\textwidth]{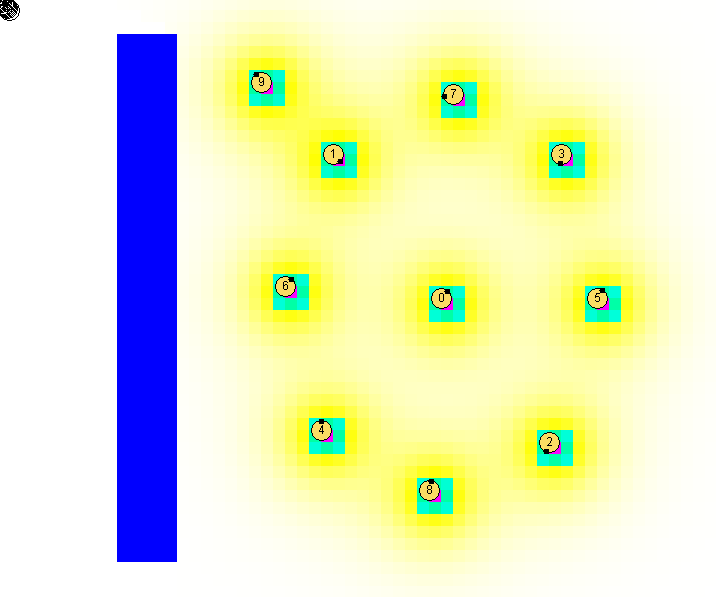}
			\caption{$H(\textbf{s})=371.9$} 	
	\end{subfigure}
\quad\begin{subfigure}[t]{1in}
			\centering
			\includegraphics[width=\textwidth]{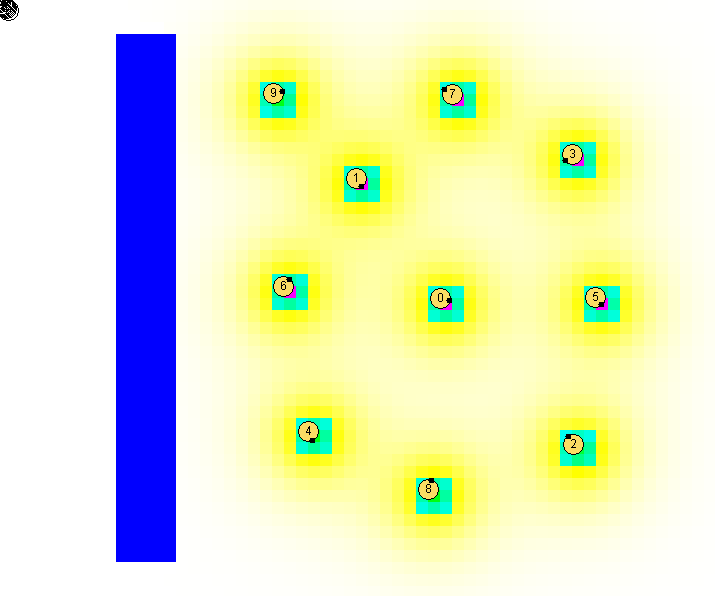}	
			\caption{$H(\textbf{s})=373.2$}
	 \end{subfigure}
\caption{The decay factor $\lambda=0.4$, and a wall-like obstacle in the
mission space}%
\label{fig:NarrowP3}%
\end{figure}

\begin{figure}[ptb]
\centering
\begin{subfigure}[t]{1in}
			\centering
			\includegraphics[width=\textwidth]{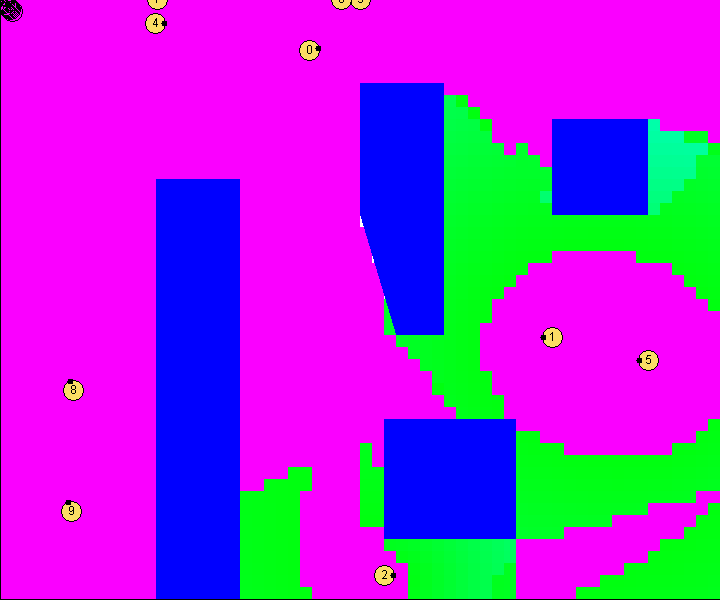}	
			\caption{$H(\textbf{s})=2401.6$}
	\end{subfigure}
\quad\begin{subfigure}[t]{1in}
			\centering
			\includegraphics[width=\textwidth]{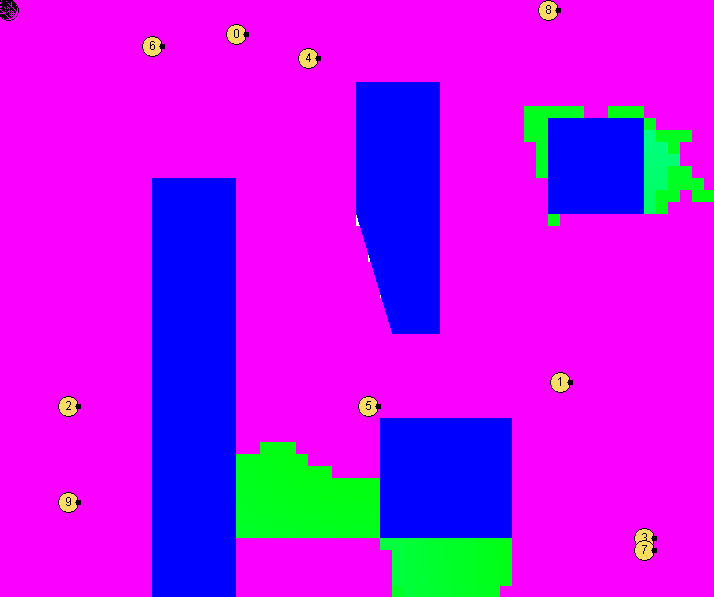}	
			\caption{$H(\textbf{s})=2421.9$} 	
	\end{subfigure}
\quad\begin{subfigure}[t]{1in}
			\centering
			\includegraphics[width=\textwidth]{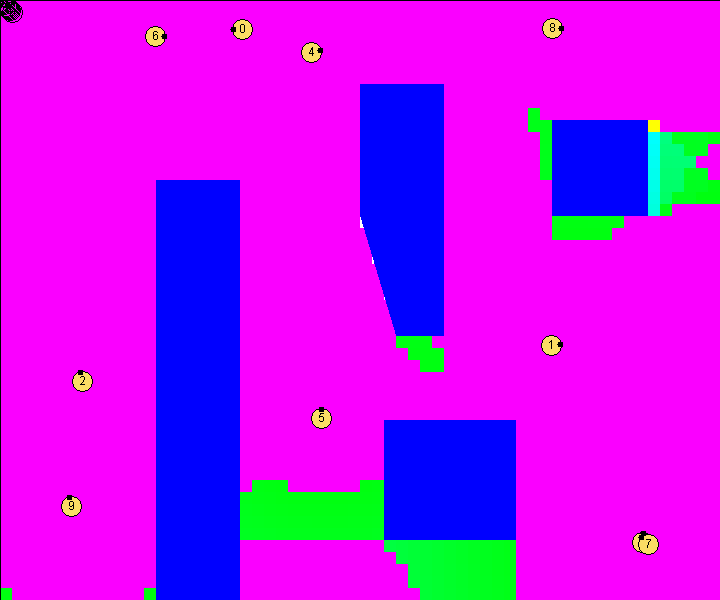}	
			\caption{$H(\textbf{s})=2423.4$}
	 \end{subfigure}
\caption{The decay factor $\lambda=0.02$, in a general mission space}%
\label{fig:GeneralP1}%
\end{figure}

\begin{figure}[ptb]
\centering
\begin{subfigure}[t]{1in}
			\centering
			\includegraphics[width=\textwidth]{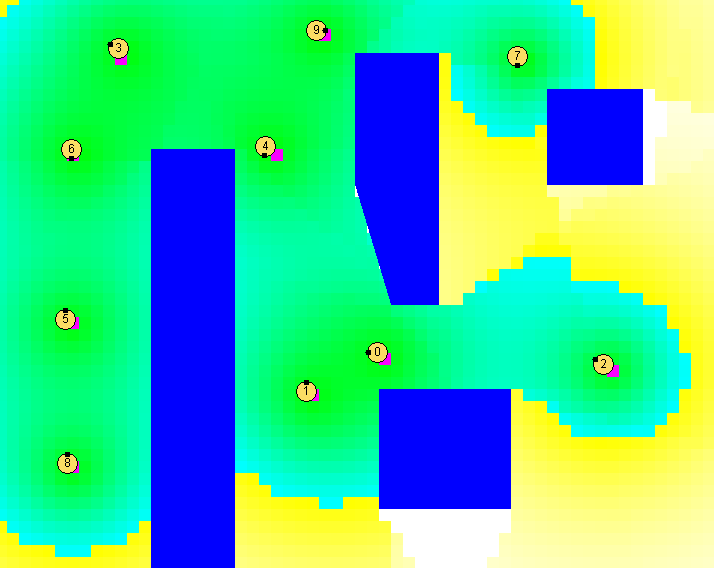}	
			\caption{$H(\textbf{s})=1443.4$}
	 \end{subfigure}
\quad\begin{subfigure}[t]{1in}
			\centering
			\includegraphics[width=\textwidth]{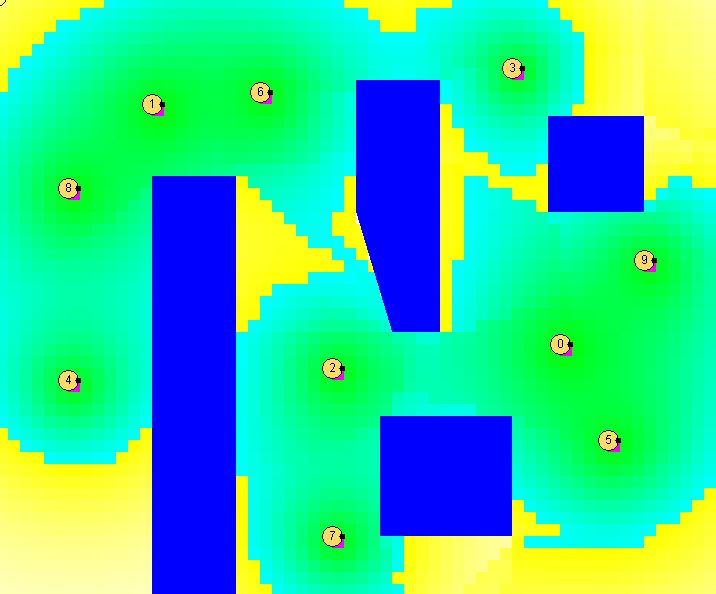}
			\caption{$H(\textbf{s})=1518.9$} 	
	\end{subfigure}
\quad\begin{subfigure}[t]{1in}
			\centering
			\includegraphics[width=\textwidth]{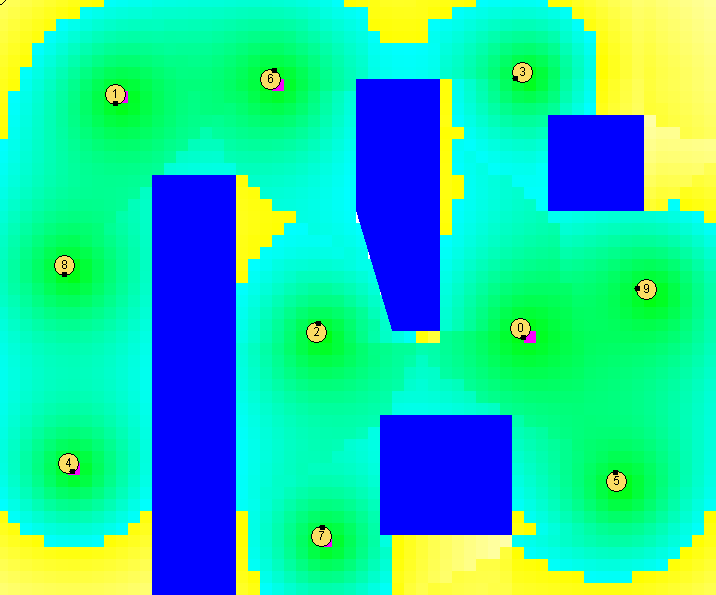}	
			\caption{$H(\textbf{s})=1532.9$}
	 \end{subfigure}
\caption{The decay factor $\lambda=0.12$, in a general mission space}%
\label{fig:GeneralP2}%
\end{figure}

\begin{figure}[ptb]
\centering
\begin{subfigure}[t]{1in}
			\centering
			\includegraphics[width=\textwidth]{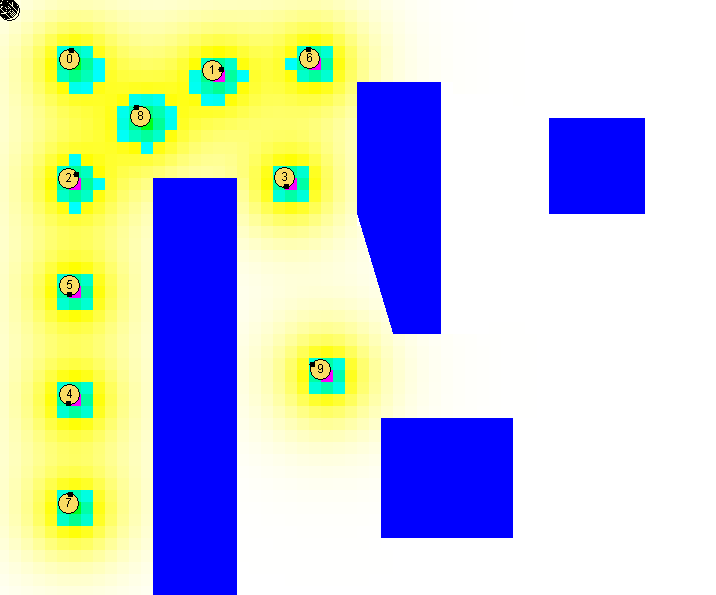}	
			\caption{$H(\textbf{s})=325.8$}
	 \end{subfigure}
\quad\begin{subfigure}[t]{1in}
			\centering
			\includegraphics[width=\textwidth]{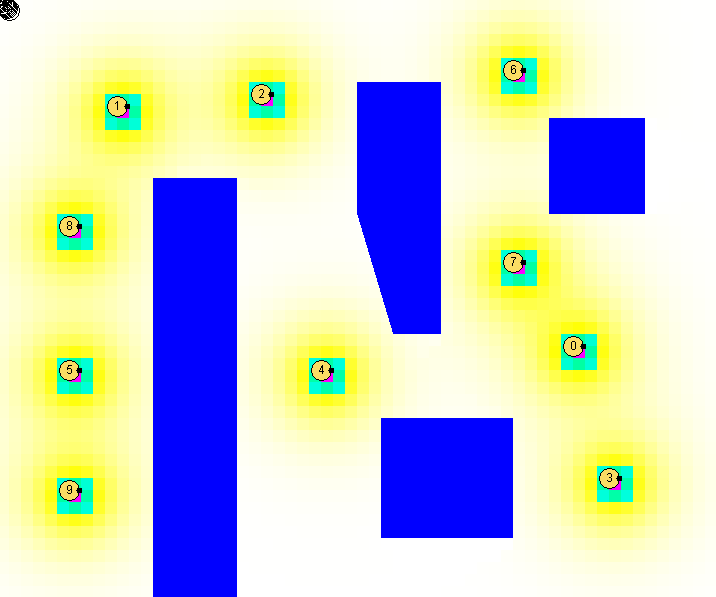}
			\caption{$H(\textbf{s})=349.2$} 	
	\end{subfigure}
\quad\begin{subfigure}[t]{1in}
			\centering
			\includegraphics[width=\textwidth]{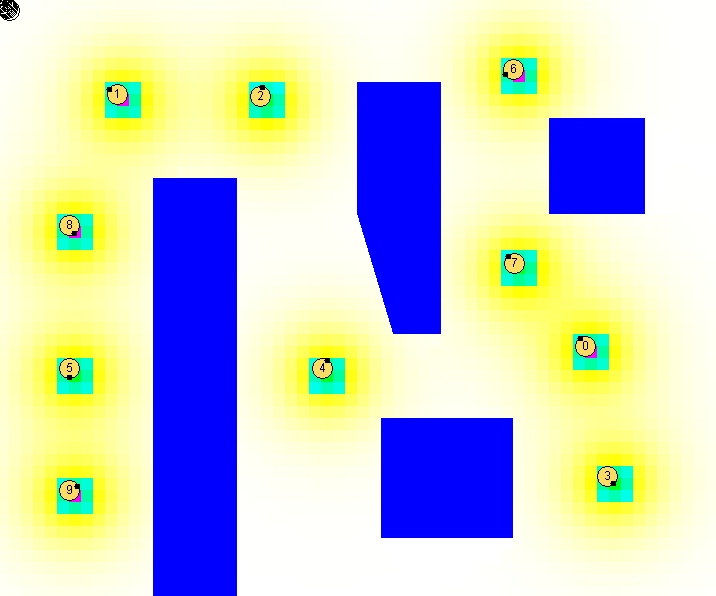}	
			\caption{$H(\textbf{s})=349.4$}
	 \end{subfigure}	
\caption{The decay factor $\lambda=0.4$, in a general mission space}%
\label{fig:GeneralP3}%
\end{figure}

\begin{figure}[ptb]
\centering
\begin{subfigure}[t]{1in}
			\centering
			\includegraphics[width=\textwidth]{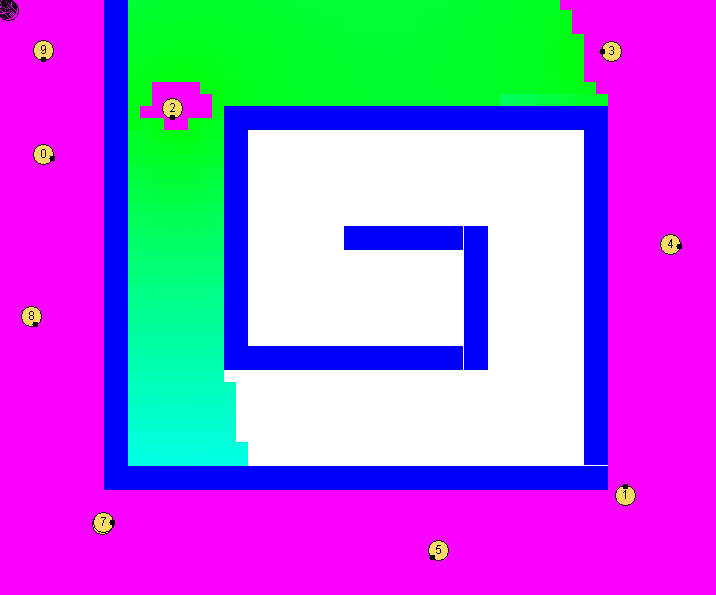}	
			\caption{$H(\textbf{s})=1792.2$}
	 \end{subfigure}
\quad\begin{subfigure}[t]{1in}
			\centering
			\includegraphics[width=\textwidth]{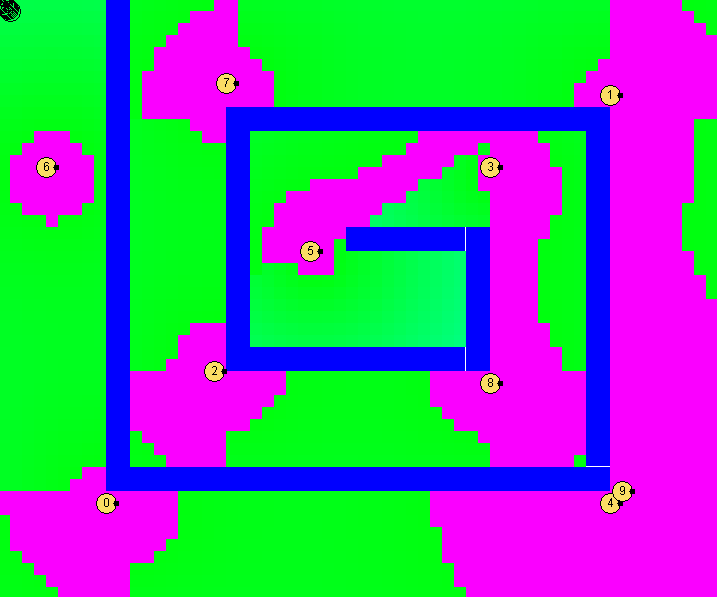}	
			\caption{$H(\textbf{s})=2490.0$} 	
	\end{subfigure}
\quad\begin{subfigure}[t]{1in}
			\centering
			\includegraphics[width=\textwidth]{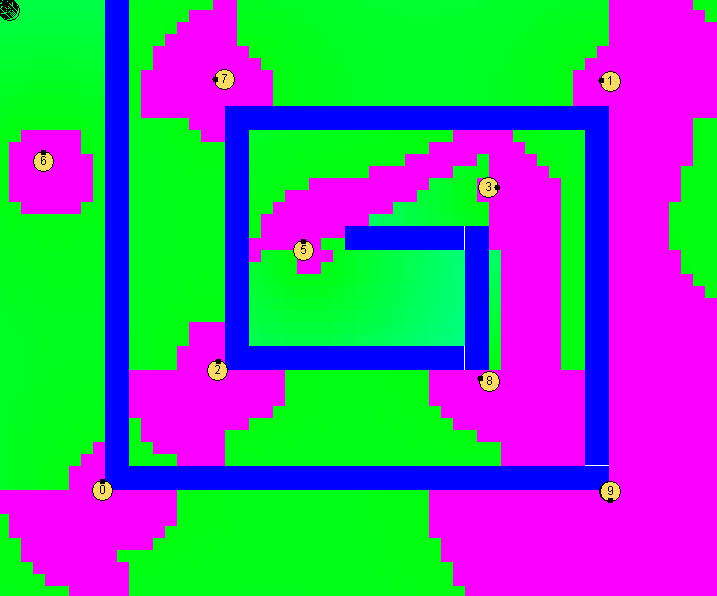}	
			\caption{$H(\textbf{s})=2490.6$}
	 \end{subfigure}
\caption{The decay factor $\lambda=0.02$, in a maze mission space}%
\label{fig:MazeP1}%
\end{figure}

\begin{figure}[ptb]
\centering
\begin{subfigure}[t]{1in}
			\centering
			\includegraphics[width=\textwidth]{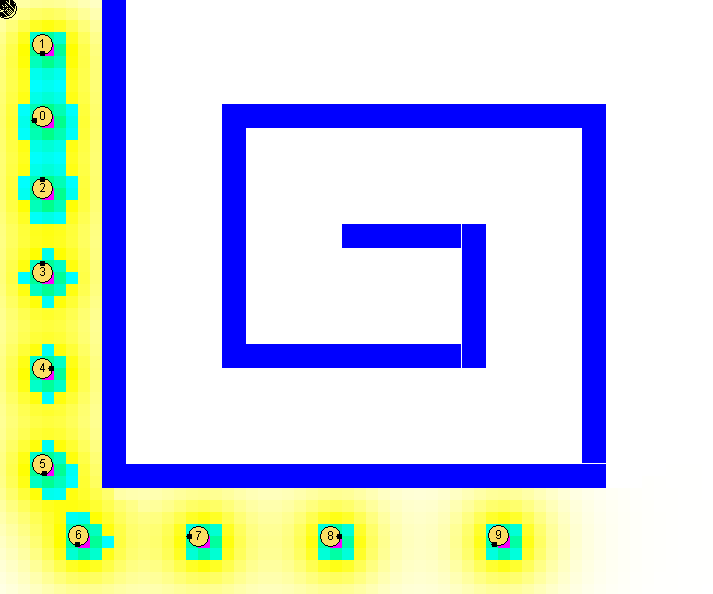}	
			\caption{$H(\textbf{s})=275.0$}
	 \end{subfigure}
\quad\begin{subfigure}[t]{1in}
			\centering
			\includegraphics[width=\textwidth]{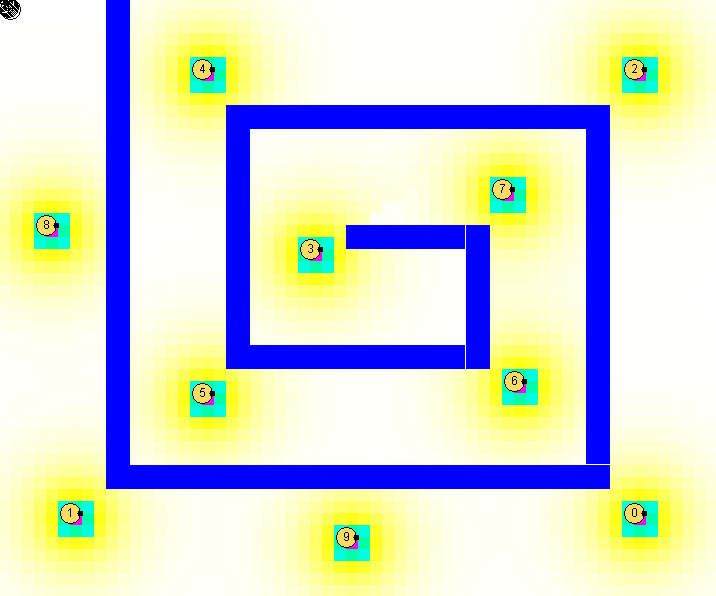}
			\caption{$H(\textbf{s})=311.1$} 	
	\end{subfigure}
\quad\begin{subfigure}[t]{1in}
			\centering
			\includegraphics[width=\textwidth]{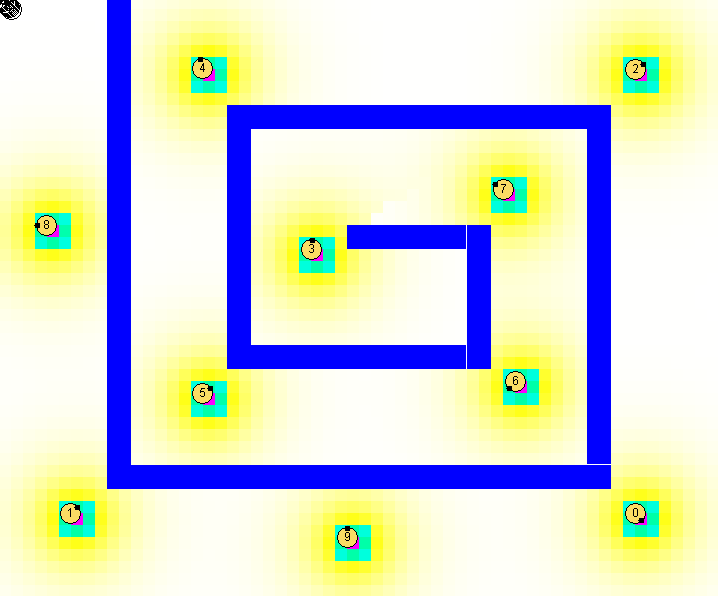}	
			\caption{$H(\textbf{s})=311.1$}
	 \end{subfigure}
\caption{The decay factor $\lambda=0.4$, in a maze mission space}%
\label{fig:MazeP3}%
\end{figure}

\begin{figure}[ptb]
\centering
\begin{subfigure}[t]{1in}
			\centering
			\includegraphics[width=\textwidth]{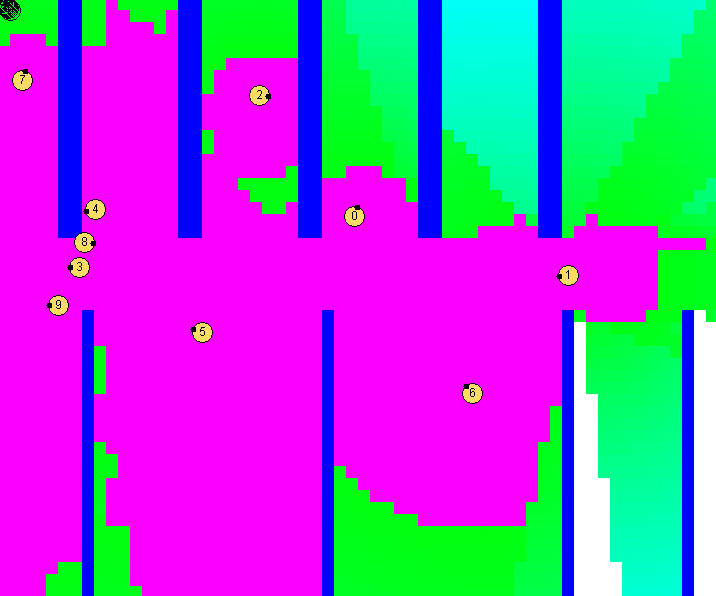}	
			\caption{$H(\textbf{s})=2418.8$}
	 \end{subfigure}
\quad\begin{subfigure}[t]{1in}
			\centering
			\includegraphics[width=\textwidth]{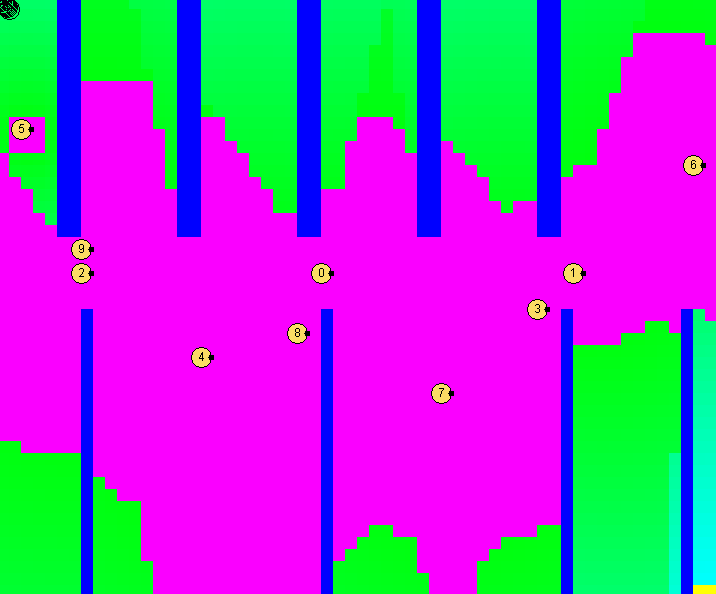}	
			\caption{$H(\textbf{s})=2582.3$} 	
	\end{subfigure}
\quad\begin{subfigure}[t]{1in}
			\centering
			\includegraphics[width=\textwidth]{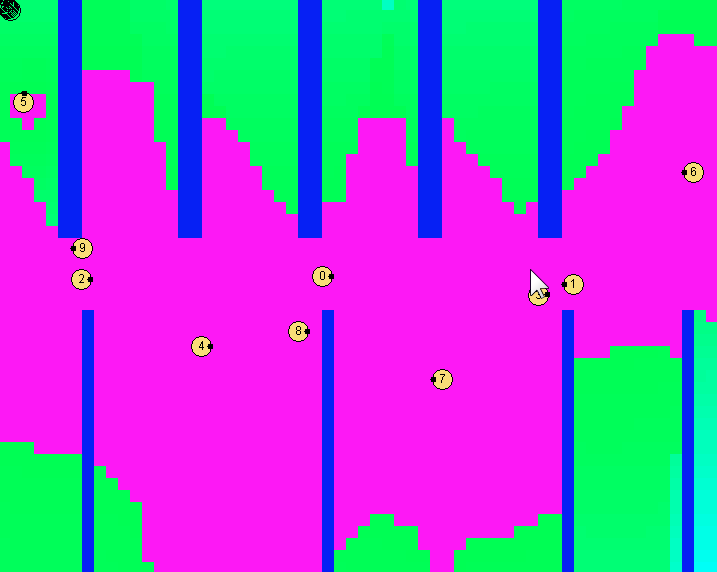}	
			\caption{$H(\textbf{s})=2583.5$}
			\label{fig:RoomGreedyGradientP1}
	 \end{subfigure}
\caption{The decay factor $\lambda=0.02$, in a room mission space}%
\label{fig:RoomP1}%
\end{figure}

\begin{figure}[ptb]
\centering
\begin{subfigure}[t]{1in}
			\centering
			\includegraphics[width=\textwidth]{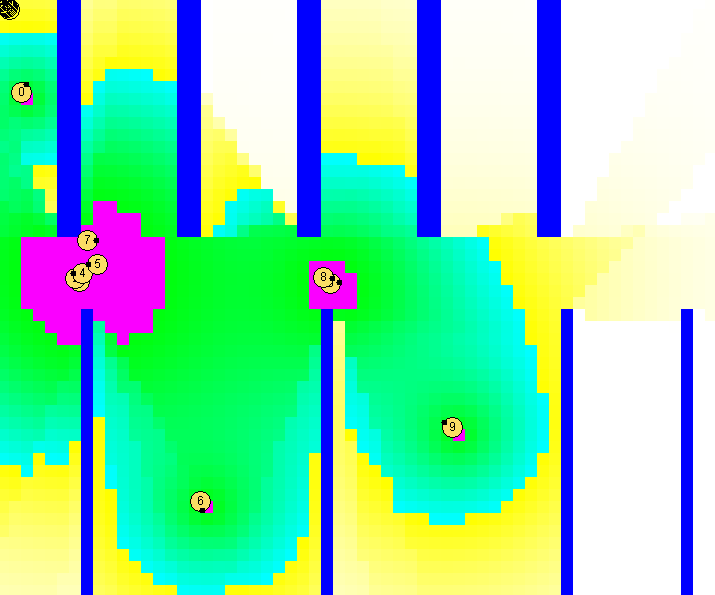}	
			\caption{$H(\textbf{s})=1187.0$}
	 \end{subfigure}
\quad\begin{subfigure}[t]{1in}
			\centering
			\includegraphics[width=\textwidth]{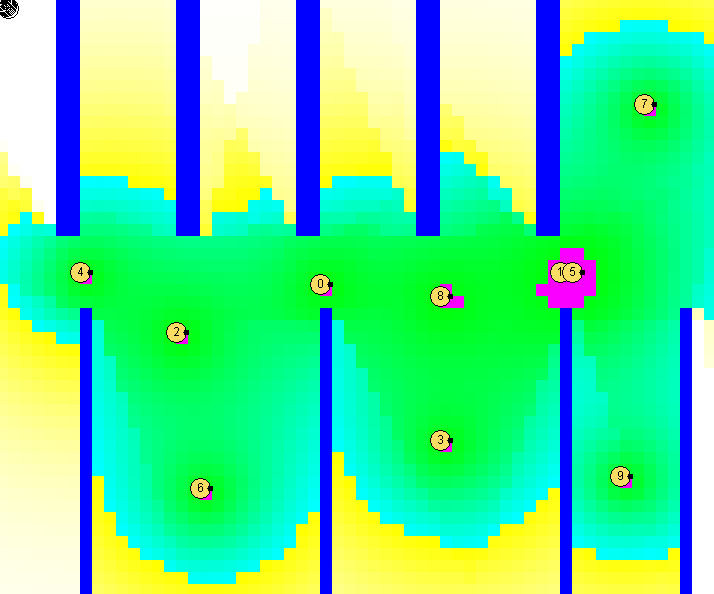}
			\caption{$H(\textbf{s})=1462.6$} 	
	\end{subfigure}
\quad\begin{subfigure}[t]{1in}
			\centering
			\includegraphics[width=\textwidth]{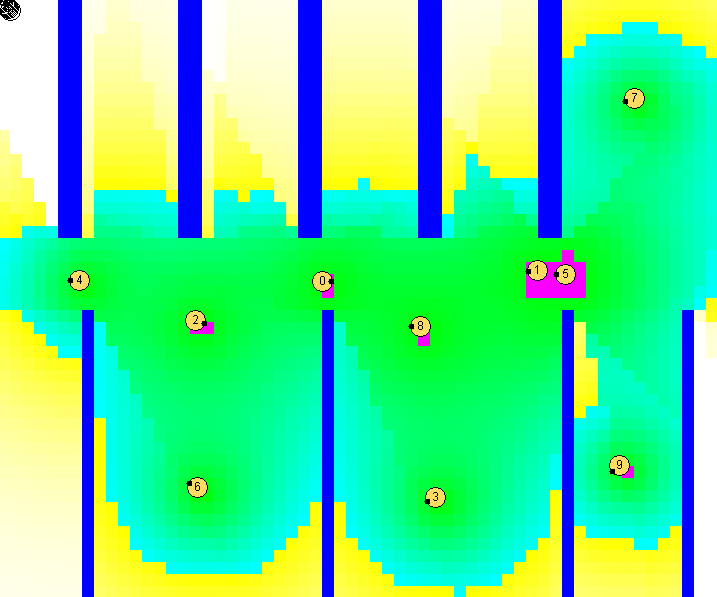}	
			\caption{$H(\textbf{s})=1466.9$}
			\label{fig:RoomGreedyGradientP2}
	 \end{subfigure}
\caption{The decay factor $\lambda=0.12$, in a room mission space}%
\label{fig:RoomP2}%
\end{figure}

\begin{figure}[ptb]
\centering
\begin{subfigure}[t]{1in}
			\centering
			\includegraphics[width=\textwidth]{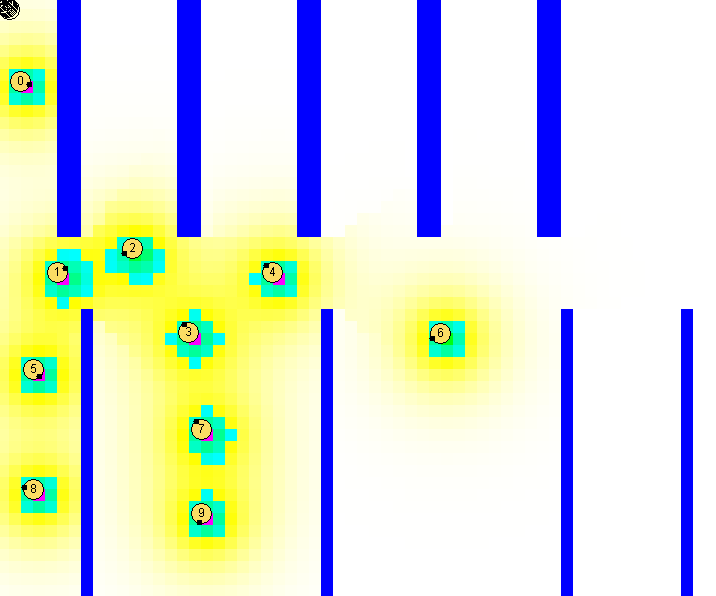}	
			\caption{$H(\textbf{s})=303.1$}
	 \end{subfigure}
\quad\begin{subfigure}[t]{1in}
			\centering
			\includegraphics[width=\textwidth]{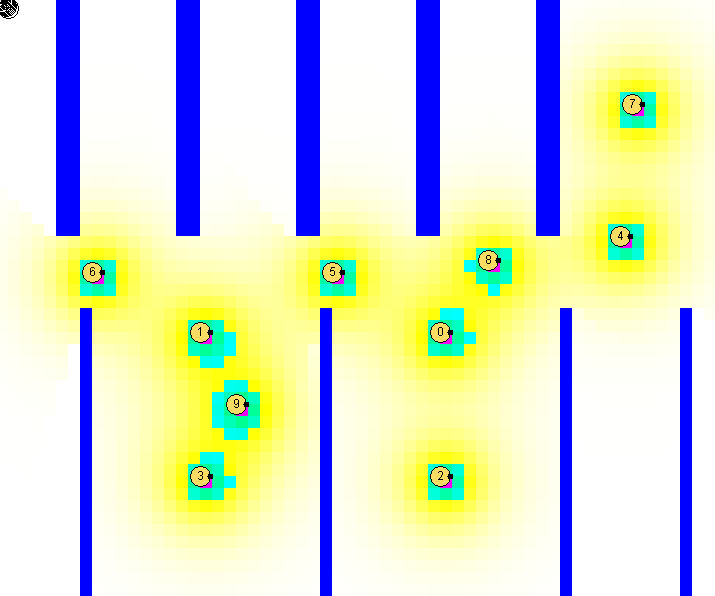}
			\caption{$H(\textbf{s})=344.5$} 	
	\end{subfigure}
\quad\begin{subfigure}[t]{1in}
			\centering
			\includegraphics[width=\textwidth]{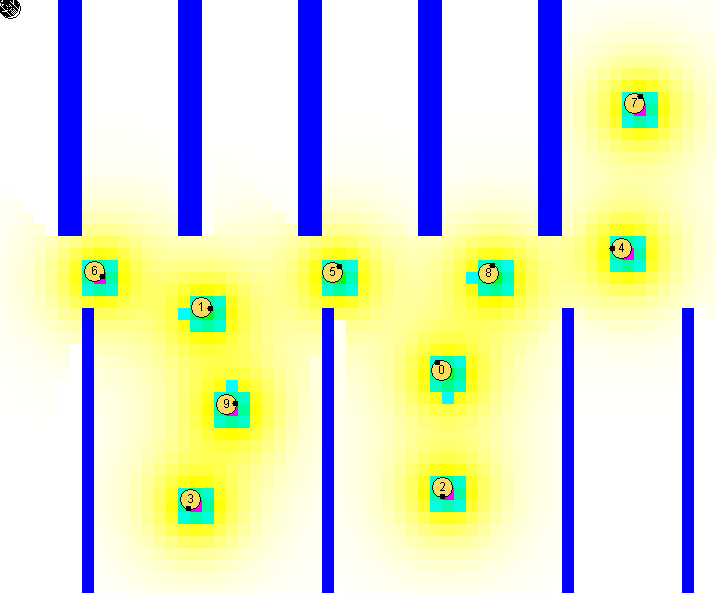}	
			\caption{$H(\textbf{s})=347.2$}
	 \end{subfigure}
\quad\caption{The decay factor $\lambda=0.4$, in a room mission space}%
\label{fig:RoomP3}%
\end{figure}

Next, we compare the performance of the greedy algorithm (Algorithm 1), the
proposed GGA (Algorithm 2) and the distributed gradient algorithm in
(\ref{gradientbasedalgo}) for solving the optimal coverage problem in
different mission spaces: no obstacles, a wall-like obstacle, a maze-like
obstacle, a collection of random obstacles, and a mission space resembling a
building with multiple rooms. Since the global optimum is unknown, we resort
to comparing all three results as shown in Figs. \ref{fig:NoObsP1}%
-\ref{fig:NoObsP3}, Figs. \ref{fig:NarrowP1}-\ref{fig:NarrowP3}, Figs.
\ref{fig:GeneralP1}-\ref{fig:GeneralP3}, Figs. \ref{fig:MazeP1}%
-\ref{fig:MazeP3} and Figs. \ref{fig:RoomP1}-\ref{fig:RoomP3} for each of
these five cases. In each case, we fix the sensing range to $\delta
_{i}=80,i=1,...,N$ and use three different values of $\lambda$, where (a)
shows the results of our distributed gradient-based algorithm, (b) shows the
results under the greedy algorithm, and (c) shows the results under the GGA.
The mission space is colored from dark to light as the joint detection
probability (our objective function) decreases: the joint detection
probability is $\geq0.97$ for purple areas, $\geq0.50$ for green areas, and
near zero for white areas.

\begin{figure}[ptb]
\centering
\begin{subfigure}[t]{1in}
			\centering
			\includegraphics[width=\textwidth]{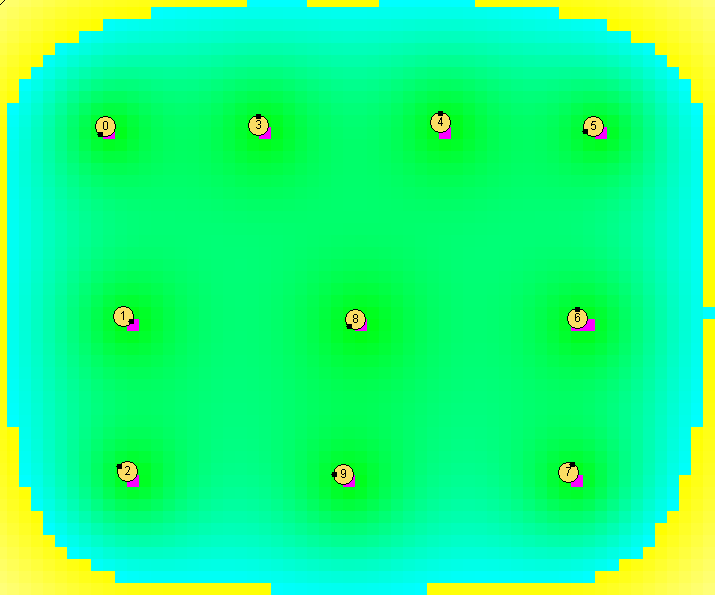}	
			\caption{$H(\textbf{s})=2105.3$}
			\label{GradientP2}
	\end{subfigure}
\quad\begin{subfigure}[t]{1in}
			\centering
			\includegraphics[width=\textwidth]{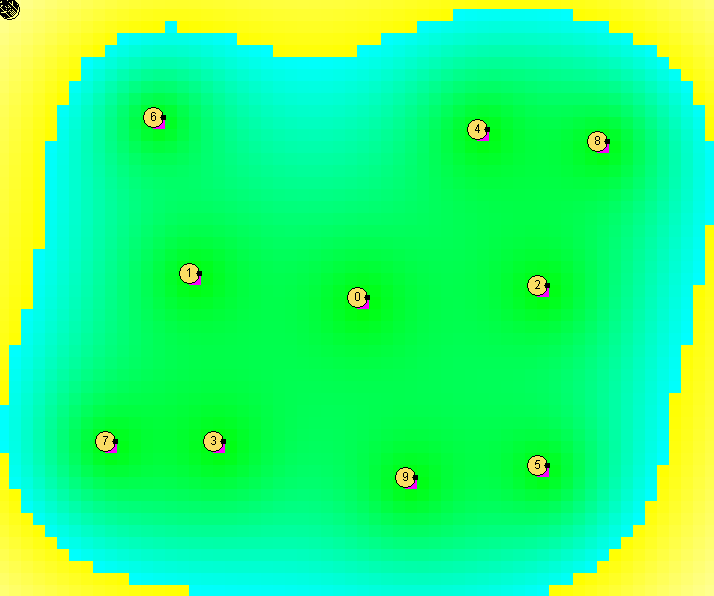}
			\caption{$H(\textbf{s})=2080.9$} 	
			\label{fig:GreedyP2}
	\end{subfigure}
\quad\begin{subfigure}[t]{1in}
			\centering
			\includegraphics[width=\textwidth]{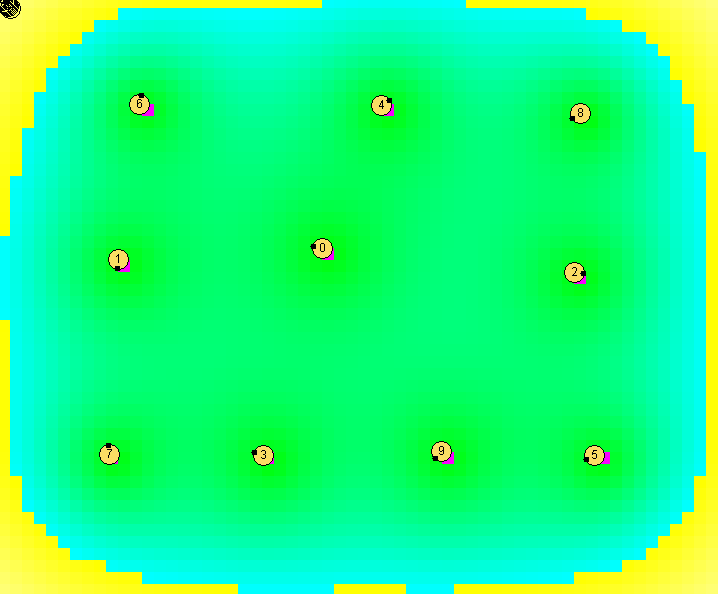}	
			\caption{$H(\textbf{s})=2105.3$}
			\label{fig:GreedyGradientP2}
	 \end{subfigure}
\quad\caption{The decay factor $\lambda=0.12$, and no obstacles in the mission
space}%
\label{fig:NoObsP2}%
\end{figure}

When there are no obstacles, all algorithms perform similarly, as shown in
Figs. \ref{fig:NoObsP1}-\ref{fig:NoObsP3}, although the actual agent
configurations are generally different (suggesting that there are multiple
equivalent local, and possibly global, optima.) For the case where
$\lambda=0.02$, the greedy algorithm is guaranteed to be within about $8\%$ of
the global optimum of (\ref{covcontrolH}) (using Fig. \ref{fig:Bound}) and we
see that using the GGA hardly improves performance, probably because the
actual global optimum is achieved.

For all cases with obstacles in the mission space, the greedy algorithm and
the GGA clearly outperform the basic gradient-based algorithm. Moreover, the
results of the GGA significantly improve upon those reported in our previous
work \cite{Sun2014}. As an example, in the cases of Figs. \ref{fig:RoomP1}%
-\ref{fig:RoomP3} with $\lambda=0.12$, the objective function value is
improved from a value of $1419.5$ reported in \cite{Sun2014} (using the
distributed gradient-based algorithm with improvements provided through the
use of boosting functions) to $1466.9$ using the GGA as shown in Fig.
\ref{fig:RoomP2}.







\begin{figure}[ptb]
\centering
\begin{subfigure}[t]{1in}
			\centering
			\includegraphics[width=\textwidth]{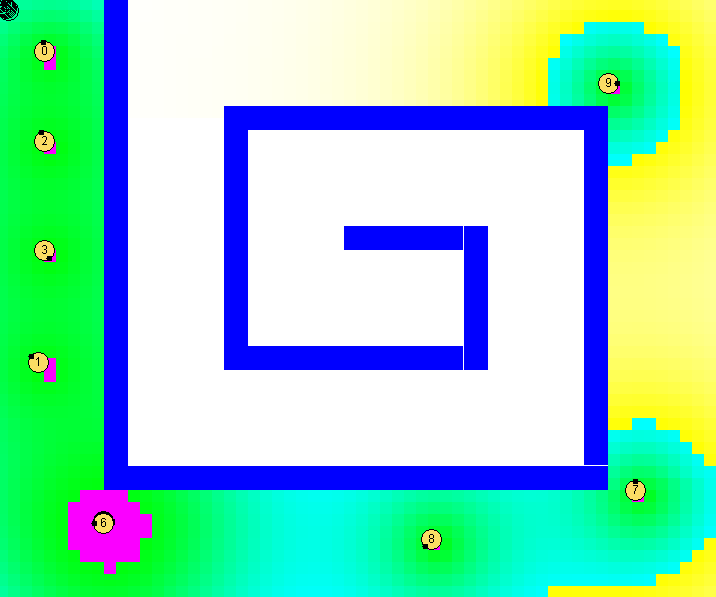}	
			\caption{$H(\textbf{s})=924.5$}
	 \end{subfigure}
\quad\begin{subfigure}[t]{1in}
			\centering
			\includegraphics[width=\textwidth]{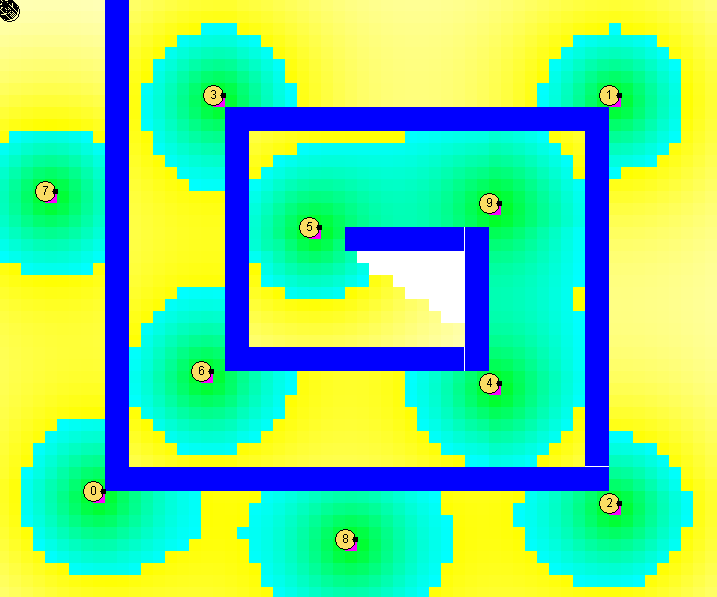}
			\caption{$H(\textbf{s})=1297.6$} 	
	\end{subfigure}
\quad\begin{subfigure}[t]{1in}
			\centering
			\includegraphics[width=\textwidth]{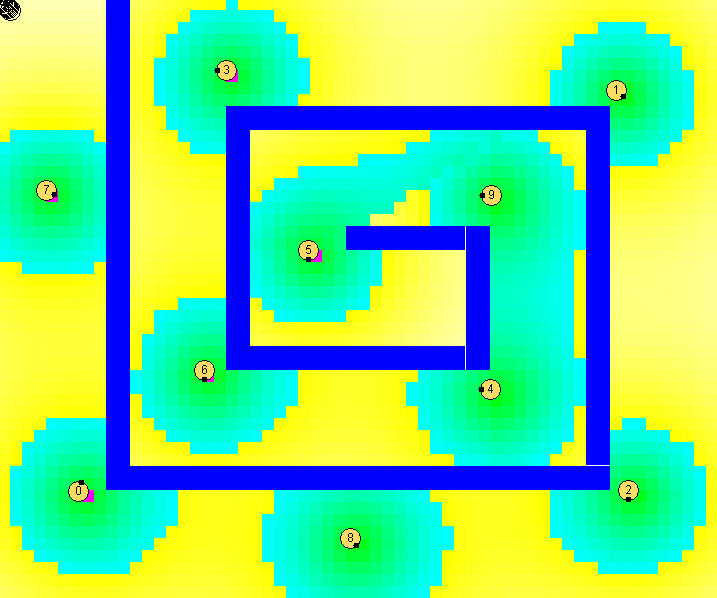}	
			\caption{$H(\textbf{s})=1307.9$}
	 \end{subfigure}
\caption{The decay factor $\lambda=0.12$, in a maze mission space}%
\label{fig:MazeP2}%
\end{figure}



\section{CONCLUSION AND FUTURE WORK}

We have obtained a solution to the optimal coverage problem through the greedy
algorithm (whose time complexity is linear in the number of agents in the
network) with a guaranteed lower bound relative to the global optimum which is
significantly tighter than the one well-known in the literature to be $1-1/e$.
This is made possible by proving that our coverage metric is monotone
submodular and by calculating its total curvature and its elemental curvature.
Therefore, we are able to reduce the theoretical performance gap between
optimal and suboptimal solutions enabled by the submodularity theory.
Moreover, we have shown that the two new bounds derived are complementary with
respect to the sensing capabilities of the agents and each one approaches its
maximal value of 1 under different conditions on the sensing capabilities,
enabling us to select the most appropriate one depending on the
characteristics of the agents at our disposal. In addition, by combining the
greedy algorithm with a distributed gradient-based algorithm we have proposed
a greedy-gradient algorithm (GGA) so as to improve the coverage performance by
searching in a continuous feasible region with initial conditions provided by
the greedy algorithm. We have included simulation results uniformly showing
that the proposed distributed GGA outperforms other related methods we are
aware of.

An interesting future research direction is to study whether a distributed
greedy algorithm can be developed and whether the lower bounds obtained
through the associated curvatures are still as tight as those we have obtained
so far.

\appendices

\section{Proof of Equivalence}

\label{app1} \textbf{Definition 1} $\Rightarrow$ \textbf{Definition 2}

Suppose that $S\subseteq T$, $y\notin T$, and $f$ satisfies (\ref{def1}).
Replacing $S$ in (\ref{def1}) by $S\cup\left\{  y\right\}  $ gives%
\begin{equation}
f(S\cup\left\{  y\right\}  \cup T)+f(S\cup\left\{  y\right\}  \cap T)\leq
f(S\cup\left\{  y\right\}  )+f(T). \label{appe1}%
\end{equation}
Rearranging the terms in (\ref{appe1}), we obtain (\ref{def2}).

\textbf{Definition 2} $\Rightarrow$ \textbf{Definition 1}

Suppose that $f$ satisfies (\ref{def2}), and $y\in S/(S\cap T)$. It is easy to
verify that $S\cap T\subseteq T$, and $y\notin T$. Replacing $S$ in
(\ref{def2}) by $S\cap T$, and using (\ref{def2}) repeatedly to all elements
$y \in S/(S\cap T)$, we obtain%
\begin{equation}%
\begin{split}
\label{appe2} &  \underset{f\left(  S\right)  }{\underbrace{f\left(  S\cap
T\cup(S/(S\cap T))\right)  }}-f\left(  S\cap T\right) \\
&  \qquad\qquad\qquad\geq\underset{f\left(  T\cup S\right)  }{\underbrace
{f\left(  T\cup(S/(S\cap T))\right)  }}-f\left(  T\right)  .
\end{split}
\end{equation}
Rearranging the terms in (\ref{appe2}) gives (\ref{def1}).
\bibliographystyle{IEEEtran}
\bibliography{research2017CDC}

\end{document}